\documentclass[12pt]{article}
\usepackage[margin=1in]{geometry}
\usepackage{amsmath,amsthm,amssymb}
\usepackage{graphicx}
\usepackage{indentfirst}

\usepackage{footnote}
\usepackage{algorithm}
\usepackage{algorithmicx}
\usepackage{algpseudocode}
\algrenewcommand\algorithmicrequire{\textbf{Input:}}
\algrenewcommand\algorithmicensure{\textbf{Output:}}

\usepackage{url}

\usepackage{float}
\usepackage{color}
\usepackage{cleveref}
\usepackage{subfig}
\usepackage{adjustbox}
\usepackage{booktabs}

\usepackage[sort,numbers]{natbib} 

\newtheorem{theorem}{Theorem}
\newtheorem{remark}[theorem]{Remark}

\crefname{assumption}{Assumption}{Assumptions}
\Crefname{assumption}{Assumption}{Assumptions}
\crefname{lemma}{Lemma}{Lemmas}
\Crefname{lemma}{Lemma}{Lemmas}
\crefname{appendix}{Appendix}{Appendices}
\Crefname{appendix}{Appendix}{Appendices}
\crefname{equation}{}{}
\Crefname{equation}{}{}
\crefname{figure}{Figure}{Figures}
\Crefname{figure}{Figure}{Figures}

\newcommand{\R}{\mathbb{R}}
\newcommand{\di}{\,\mathrm{d}}

\newcommand{\mO}{\mathcal{O}}

\newcommand{\E}[1]{\mathbb{E}\left[#1\right]}

\newcommand{\br}[1]{\left(#1\right)}

\newcommand{\abs}[1]{\left\vert#1\right\vert}

\newcommand{\bbr}[1]{\left\{#1\right\}}

\newcommand{\B}[1]{\boldsymbol{#1}}
\newcommand{\sptext}[1]{\;\; \text{#1} \;\;}

\newcounter{jointfn}  

\begin{document}
 
%
%
 
 \title{Deep Policy Iteration for High-Dimensional Mean-Field Games with Regenerative Reformulation
}

\author{
Shuixin Fang\thanks{SKLMS \& Institute of Computational Mathematics and Scientific/Engineering Computing, Academy of Mathematics and Systems Science, Chinese Academy of Sciences, Beijing 100190, China. Email: \texttt{sxfang@amss.ac.cn}, \ \texttt{tzhou@lsec.cc.ac.cn}}\setcounter{jointfn}{\value{footnote}}\ , 
Shupeng Wang\thanks{School of Control Science and Engineering, Shandong University, Jinan 250061, China. Email: \texttt{wang\_shupeng@sdu.edu.cn}}\ , Zhen Wu\thanks{School of Mathematics, Shandong University, Jinan 250100, China. Email: \texttt{wuzhen@sdu.edu.cn}}\ , 
Hui Zhang\thanks{Department of Applied Mathematics, The Hong Kong Polytechnic University, Hung Hom, Hong Kong. Email: \texttt{hui1203.zhang@polyu.edu.hk}}\ , 
Tao Zhou\footnotemark[\arabic{jointfn}]
}
\date{}
\maketitle

\begin{abstract}
This paper develops a deep policy iteration method for high-dimensional finite-horizon mean-field games (MFG). 
We reformulate the game as a regenerative problem with deterministic cycles, which allows policy evaluation (PE), policy improvement (PI), and population measure estimation to be carried out cycle by cycle. 
Within this formulation, we approximate the population measure by a particle system and update it using a one-step random mapping induced by the Euler-Maruyama discretization of the state dynamics. 
This update transports a mini-batch of particles from one cycle to the next, avoiding sequential trajectory simulation over the entire time horizon at each iteration. The PE and PI subproblems are formulated through the relation between consecutive cycles, with adversarial training used for evaluation and averaged optimization used for improvement. The resulting method is efficient and scalable in high dimensions, as it avoids the direct solution of the coupled Hamilton-Jacobi-Bellman and Fokker-Planck system, the full simulation of trajectories to estimate the population measure, the explicit computation of conditional expectations in policy evaluation, and pointwise optimization in policy improvement. Numerical experiments demonstrate that the proposed method effectively handles dimensions up to 10,000.
\\
\textbf{key words}: 
mean-field games, high-dimensional problem, deep learning, adversarial learning, reinforcement learning, policy iteration 
\end{abstract}

\section{Introduction}

Large-population systems arise in a wide range of applications in engineering, economics, finance, and the social sciences. Some
significant progress has been made in the study of large-population systems, which involve numerous interacting agents. Although an individual agent has a negligible impact at the micro level, the collective behavior of all agents is highly significant at the macro level. Mathematically, there exists a weak coupling relation among the considerable agents in their dynamics or cost functionals through the state average or more general empirical measure. However, when the number of agents becomes sufficiently large, the problem becomes intractable to deal with directly due to the curse of dimensionality and the weak coupling structure. For these reasons, it is unrealistic for an individual agent to obtain a centralized strategy, as it would require the information of other peers. A more feasible alternative is to design decentralized control strategies based solely on local information and some precomputable off-line quantities. It should be emphasized that the MFG approach plays a key role in achieving this goal. For a comprehensive understanding of MFG theory, interested readers can refer to \cite{Bensoussan13,Bensoussan18Mean,cardaliaguet2010notes,Carmona18} for a detailed overview.

MFGs describe Nash equilibria in differential games with an infinite number of players. The theory captures the collective behavior emerging from a continuum of rational agents, each seeking to minimize a cost or maximize a benefit. This cost depends on the state of the agent and also on the statistical distribution of all the agents' states, that is, the global state of the system. MFG theory has garnered significant attention, particularly due to pioneering contributions from Lasry and Lions \cite{Lasry13}, Huang et al. \cite{Minyi10,Minyi12}, Andersson and Djehiche \cite{Andersson}, Yong \cite{Yong13}, and Elliott et al. \cite{Elliott13}, fostering a significant development in MFG. It has found applications in finance and economy \cite{Carmona15probabilistic,Lachapelle16,Lachapelle10}, systemic risk \cite{Carmona15Mean,Garnier13}, crowd motion \cite{Achdou19,lachapelle2011mean}, and so on. Compared with conventional stochastic optimal control problems (SOCPs), MFGs feature that both the state dynamics and the agents' cost functionals depend on the law of the state process. In the finite-horizon setting, dynamic programming for MFG typically yields a coupled system consisting of a Hamilton-Jacobi-Bellman (HJB) equation and a Fokker-Planck (FP) equation. In high-dimensional settings, numerical methods for this coupled HJB-FP system suffer from the curse of dimensionality (CoD).

In recent years, deep learning methods have shown great promise in alleviating the CoD. A substantial body of work has emerged on deep learning approaches for solving high-dimensional PDEs, including physics-informed neural networks (PINNs) \cite{Raissi2019Physics}, deep convolution residual neural networks \cite{Ruthotto20}, deep Galerkin method \cite{Sirignano2018DGM}, deep BSDE methods \cite{han2018solving,Zhang2022FBSDE}, deep backward schemes \cite{germain2022Approximation,hure2020deep}, DeepMartNet \cite{cai2024martingale,cai2024socmartnet,Andrew26}, shot-gun method \cite{zhang2025shotgun}, deep random difference method \cite{Cai2026Deep}, etc. However, directly treating the coupled HJB-FP system remains challenging.

The main difficulty is twofold. First, the coupling between the HJB and FP equations requires repeatedly computing high-dimensional integrals with respect to the state distribution, which leads to substantial computational cost. Second, spatial sampling plays a critical role in the accuracy of deep-learning methods for PDEs and SOCPs, as demonstrated in \cite{Gao2023Failure,Tang2023DAS,Xiaoliang}. An effective sampling strategy should capture the intrinsic dynamics of the underlying problem \cite{Li2024Neural,Cai2026Deep}. However, for the coupled HJB-FP system, spatial sampling is simultaneously required to approximate the high-dimensional integrals with respect to the state distribution, enforce the HJB equation, and resolve the FP equation. Such difficulties persist even in recent machine-learning approaches for high-dimensional mean field games and mean field control \cite{Ruthotto2020Machine}. Designing a single sampling strategy that captures all three aspects is therefore highly nontrivial.

To address these difficulties, we reformulate the finite-horizon MFG as a regenerative problem with deterministic cycles and a resetting mechanism. This reformulation replaces the original finite-horizon equilibrium condition by an equivalent characterization in terms of a cycle-wise occupation measure and leads to a policy-iteration procedure in which PE, PI, and population measure estimation are all performed on a single-cycle basis. The main advantage of this viewpoint is that it avoids working directly with the coupled HJB--FP system at each iteration. Instead, the equilibrium distribution is represented through repeated short-horizon transitions, which provides a more convenient basis for sample-based approximation in high dimensions. In this way, the reformulation converts the original problem into a structure that is more compatible with iterative learning algorithms and large-scale numerical implementation.

Based on the regenerative formulation, we develop a deep policy iteration method for high-dimensional MFG. The equilibrium occupation measure is approximated by an empirical particle measure and updated by a one-step random mapping induced by the Euler--Maruyama discretization. At each iteration, only a random mini-batch of particles is transported from one cycle to the next, avoiding full trajectory simulation over the entire horizon. Moreover, the discretized PE and PI steps are formulated in weak form: the value function and feedback control are represented by neural networks, while an adversarial test-function network yields a Galerkin-type formulation for PE. This avoids direct solution of the coupled HJB--FP system, explicit computation of conditional expectations in PE, higher-order derivatives of the HJB equation, and pointwise Hamiltonian minimization in PI. In addition, the same sample transitions are used for PE, PI, and empirical measure updates, so that policy iteration, measure approximation, and network training share a common set of state transitions, further reducing the overall computational cost. Numerical results demonstrate good performance up to dimension $10{,}000$.

The rest of this paper is organized as follows. In \cref{sec_problem_ref}, we formulate the MFG problem on a finite horizon, reformulate it as a regenerative MFG, and introduce the corresponding policy iteration method.
\Cref{sec_meth} develops the deep policy iteration algorithm, including the empirical approximation of the equilibrium occupation measure, mini-batch updates based on the random mapping, and the weak-form PE/PI training procedure for the neural networks. Some numerical experiments are shown in \cref{sec_num} to support our method. Finally, the conclusion is given in \cref{sec_conclusion}.

\section{Problem setup and reformulation}\label{sec_problem_ref}

In this section, we first introduce the time-dependent MFG model and then reformulate it as a regenerative MFG model, with deterministic cycles and a resetting mechanism. The latter formulation is more convenient for the development of the numerical method in \Cref{sec_meth}.
Then we introduce the policy iteration method for solving the regenerative MFG model, which forms the basis of the numerical method proposed in \Cref{sec_meth}.

\subsection{Mean-field games in finite-horizon formulation}\label{sec_problem}

Let $(\Omega, \mathcal{F}, \mathbb{P})$ be a complete probability space equipped with a filtration $(\mathcal{F}_t)_{t \geq 0}$ generated by a standard $q$-dimensional Brownian motion $B: [0, +\infty) \times \Omega \to \R^q$, augmented by the $\sigma$-algebra generated by an initial state $X_0: \Omega \to \{0\} \times \R^d$ of a representative agent.
Denote by $\mathcal{P}([0, T] \times \R^d)$ the space of probability measures on $[0, T] \times \R^d$.
The initial state $X_0$ follows a prescribed distribution $\mu_{X_0} \in \mathcal{P}(\{0\} \times \R^d)$, and is independent of $B$.  
Here, $X_0$ and $\mu_{X_0}$ are supported on $\{0\} \times \R^d$, instead of $\R^{1+d}$, because the first component of the state process is reserved for the time variable.

We consider a continuum of identical agents whose interactions are determined by a population distribution $\mu \in \mathcal{P}([0, T) \times \R^d)$. 
The state process of a representative agent, denoted by $X^{\mu, u}: [0, T) \times \Omega \to [0, T) \times \R^d$, is governed by the mean-field stochastic differential equation (MSDE)
\begin{equation}\label{eq_state}
    X_t^{\mu, u} = X_0 + \int_0^t b\big( X_s^{\mu, u}, \mu, u(X_s^{\mu, u}) \big) \di s + \int_0^t \sigma\big( X_s^{\mu, u}, \mu, u(X_s^{\mu, u}) \big) \di B_s, \quad t \in [0, T). 
\end{equation}
In the MSDE~\eqref{eq_state}, $b$ and $\sigma$ denote the drift and diffusion coefficients, respectively. 
They are chosen such that for $t \in [0, T)$, the first component of $X_t^{\mu, u}$ acts as the time variable, namely, for $x = (t, z) \in [0, T) \times \R^d$,
\begin{equation}\label{eq_bzsgmz}
    b(x, \mu, u) = \big(1,\, b_z^{\top}(x, \mu, u)\big)^{\top}, \quad \sigma(x, \mu, u) = \big(0_q, \,\sigma_z^{\top}(x, \mu, u)\big)^{\top},
\end{equation}
where the functions $b_z$ and $\sigma_z$ take values in $\R^d$ and $\R^{d \times q}$, respectively, and $0_q$ denotes the $q$-dimensional zero vector.
The function $u$ is a feedback control belonging to the admissible control set defined as
\begin{equation}\label{eq1_defUad}
    \mathbb{U}_{\mathrm{ad}}^{\mu} := \big\{u: [0, T) \times \R^d \to U \;\big|\; \text{the MSDE~\eqref{eq_state} admits a unique strong solution } X^{\mu, u}\big\},
\end{equation}
where $U \subset \R^m$ is the control domain.
Here we incorporate the well-posedness of \eqref{eq_state} into the definition of $\mathbb{U}_{\mathrm{ad}}^{\mu}$, thereby avoiding technical complications.

The representative agent aims to choose a control $u \in \mathbb{U}_{\mathrm{ad}}^{\mu}$ that minimizes the cost functional $J(u, \mu)$, defined by
\begin{equation}\label{eq_cost}
    J(u, \mu) := \E{\int_0^T f\big( X_s^{\mu, u}, \mu, u(X_s^{\mu, u}) \big) \di s + g\big( X_{T-}^{\mu, u}, \mu \big)}
\end{equation}
where the deterministic functions $f$ and $g$, both taking values in $\R$, represent the running and terminal costs, respectively.
The terminal cost is evaluated at the left limit $X_{T-}^{\mu, u}$ of $X_t^{\mu, u}$ as $t \uparrow T$, rather than at $X_T^{\mu, u}$, 
because in the regenerative formulation introduced in the next subsection, the state process $X^{\mu, u}$ undergoes a jump at time $T$.

A mean-field game equilibrium is defined as a pair $(u^*, \mu^*) \in \mathbb{U}_{\mathrm{ad}}^{\mu^*} \times \mathcal{P}([0, T) \times \R^d)$, 
such that $u^*$ minimizes the individual cost functional and $\mu^*$ coincides with the occupation measure of the optimal state process $X^{\mu^*, u^*}$, i.e.,
\begin{equation}\label{eq_equilibrium}
    J(u^*, \mu^*) = \inf_{u \in \mathbb{U}_{\mathrm{ad}}^{\mu^*}} J(u, \mu^*), \quad \mu^*(A) = \frac{1}{T} \int_0^T \E{I(X_t^{\mu^*, u^*} \in A)} \di t, \;\; \forall A \in \mathcal{B}([0, T) \times \R^d),
\end{equation}
where $\mathcal{B}([0, T) \times \R^d)$ denotes the Borel $\sigma$-algebra on $[0, T) \times \R^d$; $I(\cdot)$ is the indicator function.  
The recovery of the optimal control $u^*$ constitutes the main focus of this work.

\begin{remark}[Connection to standard MFG formulations]
    The MFG~\eqref{eq_equilibrium} is slightly more general than the standard formulation.
    In the standard setting, the functions $b$, $\sigma$, $f$, and $g$ depend on the population distribution $\mu$ only through its time-$t$ marginal $\mu_t$ on $\R^d$; for example, one writes $b(t, z, \mu_t, u)$, and similarly for $\sigma$, $f$, and $g$. 
    By contrast, in~\eqref{eq_equilibrium}, these functions may depend on the entire distribution $\mu$ on $[0, T) \times \R^d$, thereby allowing for more general interactions among agents.
    This generalization is introduced to facilitate the regenerative formulation and the numerical method developed in later sections.
\end{remark}

\subsection{Mean-field games in regenerative formulation}\label{sec_mfg_reg}

To streamline the presentation, we recast the finite-horizon MFG~\eqref{eq_state}--\eqref{eq_equilibrium} on $[0,T)$ as an infinite-horizon regenerative process on $[0,\infty)$ by repeating length-$T$ cycles.
Let $i$ be the cycle index. 
Assume $(\mathcal{F}_t)_{t\ge 0}$ is rich enough to support a sequence of cycle-start states $\{\xi_i\}_{i=0}^{\infty}$ that are independent and identically distributed (i.i.d.) across $i$ with common distribution $\mu_{X_0}$, and independent of $B$.
For each $i = 0,1,2,\dots$, the state process $X^{\mu, u}$, within the $i$-th cycle, is defined by the following MSDE:
\begin{equation}\label{eq1_state}
    X_t^{\mu, u} = X_{i T}^{\mu, u} + \int_{i T}^t b\big( X_s^{\mu, u}, \mu, u(X_s^{\mu, u}) \big) \di s + \int_{i T}^t \sigma\big( X_s^{\mu, u}, \mu, u(X_s^{\mu, u}) \big) \di B_s,
\end{equation}
for $t \in [iT, (i+1)T)$, where a resetting mechanism is applied at the end of each cycle:
\begin{equation}\label{eq_reset}
    X_{(i+1)T}^{\mu, u} = \xi_{i+1}.
\end{equation}

By construction from \eqref{eq1_state} and \eqref{eq_reset}, the state process $X^{\mu, u}$ has the same distribution over each interval $[iT, (i+1)T)$, $i = 0, 1, 2, \cdots$.
Consequently, the cost functional in \eqref{eq_cost} and the mean-field equilibrium condition in \eqref{eq_equilibrium} can be reformulated, respectively, as a long-run average cost and an invariant occupation measure:
\begin{align}
    J(u, \mu) &= \lim_{t \to \infty} \frac{T}{t} \E{\int_{0}^{t} f\big( X_s^{\mu, u}, \mu, u(X_s^{\mu, u}) \big) \di s + \sum_{i=1}^{\infty} g\big( X_{i T-}^{\mu, u}, \mu \big) I(i T \leq t) }, \label{eq_avg_J} \\
    \mu^*(A) &= \lim_{t \to \infty} \frac{1}{t} \int_0^t \E{I(X_s^{\mu^*, u^*} \in A)} \di s, \;\; \forall A \in \mathcal{B}([0, T) \times \R^d), \label{eq_avg_mu}
\end{align}
where $I(\cdot)$ denotes the indicator function.
This regenerative formulation is closely related to policy iteration methods for solving MFGs: the cycle index $i$ in \eqref{eq1_state} and \eqref{eq_reset} plays a role analogous to the iteration index, as discussed in \cref{sec_dpp}.

\subsection{Policy iteration}\label{sec_dpp}

In this subsection, we introduce policy iteration methods for the regenerative MFG model given in \cref{sec_mfg_reg}. 
At this stage, we treat the equilibrium distribution $\mu^*$ as a given parameter, thereby reducing the problem to a standard stochastic optimal control problem.
The resulting policy iteration method is closely related to the dynamic programming principle in stochastic control and reinforcement learning; see, e.g., \cite[Chapter 4]{Yong1999Stochastic}, \cite[Chapter 3]{Pham2009Continuous}, and \cite[section 4.3]{Sutton2018Reinforcement} for foundational concepts, and \cite{Yanwei2022Policy,Yanwei2023Policy,Zhou2024Solving,Zhou2025Policy,Zhou2021Actor,Mou2024Bellman,Zhu2025Optimal} for recent developments.

The cost functional in \eqref{eq_cost} motivates the definition of the value function $v: [0, T) \times \R^d \to \R$ as the optimal cost incurred during a single cycle, i.e., for $x = (t, z) \in [0, T) \times \R^d$, 
\begin{equation}\label{eq_vx}
    v(x) := \inf_{u \in \mathbb{U}_{\mathrm{ad}}^{\mu^*}} \E{\int_t^T f\big( X_s^{\mu^*, u}, \mu^*, u(X_s^{\mu^*, u}) \big) \di s + g\big( X_{T-}^{\mu^*, u}, \mu^* \big) \,\Big|\, X_t^{\mu^*, u} = x},
\end{equation}
where $\mu^*$ is the equilibrium distribution in \eqref{eq_equilibrium}.
Given the equilibrium distribution $\mu^*$, the dynamic programming principle for minimizing $J(u,\mu^*)$ in \eqref{eq_cost} consists of the following policy-evaluation (PE) and policy-improvement (PI) steps: for $i = 0, 1, 2, \cdots$, 
\begin{alignat}{4}
    &\text{find } v_{i+1} \in \mathbb{V}^{\mu^*} &\text{ s.t. }& \mathcal{R}_h(x; u_i, v_i, \mu^*) = 0, \;\;& \forall x \in [0, T-h] \times \R^d, \label{eq_dpp_v} \\
    &\text{find } u_{i+1} \in \mathbb{U}_{\mathrm{ad}}^{\mu^*} &\text{ s.t. }& \min_{u \in \mathbb{U}_{\mathrm{ad}}^{\mu^*}} \mathcal{R}_h(x; u, v_{i+1}, \mu^*), \;\;& \forall x \in [0, T-h] \times \R^d. \label{eq_dpp_u}
\end{alignat}
Here the function $\mathcal{R}_h$ denotes the dynamic programming residual, defined as follows:
\begin{equation}\label{eq_defMh}
\begin{aligned}
    &\mathcal{R}_h(x; u, v, \mu) \\
    :=\;& \E{\int_{i T + t}^{i T + t + h} f\big( X_s^{\mu, u}, \mu, u(X_s^{\mu, u}) \big) \di s + v\big(X_{(i T + t + h)-}^{\mu, u}\big) - v(x)\,\Big|\, X_{i T + t}^{\mu, u} = x}\\
    =\;& \E{\int_{t}^{t+h} f\big( X_s^{\mu, u}, \mu, u(X_s^{\mu, u}) \big) \di s + v\big(X_{(t+h)-}^{\mu, u}\big) - v(x)\,\Big|\, X_{t}^{\mu, u} = x}
\end{aligned}
\end{equation}
for any $x = (t, z) \in [0, T-h] \times \R^d$; the second equality in \eqref{eq_defMh} follows from the identity of the distribution of $X^{\mu, u}$ across cycles.
The set $\mathbb{V}^{\mu^*}$ consists of candidate value functions satisfying the terminal condition, i.e.,
\begin{equation}\label{eq_Vmu}
    \mathbb{V}^{\mu^*} := \bbr{v \in C([0, T] \times \R^d): v(x) = g(x, \mu^*), \; x \in \{T\} \times \R^d},
\end{equation}
where $C([0, T] \times \R^d)$ denotes the space of continuous functions on $[0, T] \times \R^d$.
The sequences $(v_i)_{i \ge 0}$ and $(u_i)_{i \ge 0}$ are expected to converge, as $i \to \infty$, to the value function $v$ in \eqref{eq_vx} and the optimal control $u^*$ in \eqref{eq_equilibrium}, respectively.

\section{Deep policy iteration}\label{sec_meth}

Implementing the policy iteration method in \cref{sec_dpp} requires several discretization ingredients: parameterized function classes for $v_i$ and $u_i$, an empirical approximation of the state distribution $\mu^*$, a discretization of the state space $[0, T) \times \R^d$ to enforce \eqref{eq_dpp_v} and \eqref{eq_dpp_u}, and a numerical approximation of the expectation in \eqref{eq_defMh}.
In this section, we address these issues and develop a practical algorithm for solving the MFG equilibrium problem \eqref{eq_equilibrium}.

To approximate the optimal feedback control $u^*$ in \eqref{eq_dpp_u} and the value function $v$ in \eqref{eq_vx}, we employ two neural networks, $u_{\alpha}$ and $v_{\theta}$, parameterized by $\alpha$ and $\theta$, respectively. 
The architectures of $u_{\alpha}$ and $v_{\theta}$ are described in \cref{sec_net}.
We also introduce a uniform partition of $[0, T)$ with step size $h > 0$:
\begin{equation}\label{eq_Pih}
    \Pi_h := \{n h: n = 0, 1, 2, \cdots, N-1\}, \quad \bar{\Pi}_h := \Pi_h \cup \{T\}, \quad N := T / h. 
\end{equation}
The terminal time $T$ is excluded from $\Pi_h$ and treated separately because of the resetting mechanism~\eqref{eq_reset}.

\subsection{State dynamics approximated by Markovian transition}

We first introduce a Markovian transition that describes the one-step evolution of the state process $X^{\mu, u}$ in \eqref{eq1_state} and \eqref{eq_reset}.
Let $\Phi$ be a random mapping on $\bar{\Pi}_h \times \R^d$, defined by
\begin{equation}\label{eq_randmap}
    \Phi(x, \mu, u) \overset{\mathrm{L}}{=} \begin{cases}
    x + b\big(x, \mu, u(x) \big) h + \sigma\big(x, \mu, u(x) \big) \sqrt{h} \, \zeta, \quad x \in \Pi_h \times \R^d,\\[1ex]
    X_0, \quad x \in \{T\} \times \R^d, \end{cases}
\end{equation}
where $\zeta$ is a standard Gaussian random vector in $\R^q$.
In \eqref{eq_randmap}, the notation $\overset{\mathrm{L}}{=}$ denotes equality in law; that is, $\Phi(x, \mu, u)$ is a random variable distributed as the expression on the right-hand side.

The state dynamics can be approximated by iterating the random mapping $\Phi$ over a particle system.
Specifically, let $\hat{\mu}_0$ be an initial guess of the equilibrium measure $\mu^*$ in \eqref{eq_equilibrium}, with support on $\Pi_h \times \R^d$.
Let $\{X_0^m\}_{m=1}^M$ be i.i.d. samples drawn from $\hat{\mu}_0$, yielding the initial empirical measure $\mu_0 := 1/ M \sum_{m=1}^M \delta_{X_0^m}$, where $\delta_x$ denotes the Dirac measure at $x \in \R^d$.
Let $i$ be the iteration index. 
For $i = 1, 2, \cdots$, randomly select a mini-batch of indices $\mathbb{A}_i \subset \{1, 2, \cdots, M\}$ and update the samples $\{X_i^m\}_{m=1}^M$ by
\begin{equation}\label{eq_Xnm}
    X_i^m := \begin{cases}
        \Phi(X_{i-1}^m, \mu_{i-1}, u_{\alpha}), & m \in \mathbb{A}_i,\\[1ex]
        X_{i-1}^m, & m \in \mathbb{A}_i^{c},
    \end{cases}
\end{equation}
Then define the corresponding empirical measure by
\begin{equation}\label{eq_mun}
    \mu_i := \frac{1}{M} \sum_{m=1}^M \delta_{X_i^m}. 
\end{equation}
For rigor, we assume that the random mappings $\Phi$ in \eqref{eq_randmap} are independent across particles $m=1, 2, \cdots$ and time steps $i=1,2,\cdots$, and are also independent of $\{X_0^m\}_{m=1}^M$.
Under this construction, the sample ensemble $\{X_i^m\}_{m=1}^M$ forms a discrete-time Markov chain in the time index $i$.

Combining \eqref{eq_Xnm} and \eqref{eq_mun} yields an update rule for both the state samples $\{X_i^m\}_{m=1}^M$ and the empirical measure $\mu_i$ from iteration $i-1$ to $i$.
The PE and PI steps in \eqref{eq_dpp_v} and \eqref{eq_dpp_u} can then be carried out with the equilibrium distribution $\mu^*$ replaced by the empirical measure $\mu_i$ at iteration $i$.
As the iterations proceed, the control network $u_{\alpha}$ is expected to provide a good approximation of the optimal control $u^*$, while the empirical measure $\mu_i$ is expected to converge to the equilibrium measure $\mu^*$ in \eqref{eq_equilibrium} as $i, M \to \infty$.
This forms the numerically tractable counterparts of \eqref{eq_dpp_v} and \eqref{eq_dpp_u}, as described in the next subsection.

\begin{remark}[Consistency of the random mapping]
    The random mapping $\Phi$ restricted to $x \in \Pi_h \times \R^d$ can be viewed as a one-step Euler--Maruyama approximation of the MSDE~\eqref{eq1_state}. 
    Under sufficient regularity conditions, it attains a local truncation error of order $\mO(h^2)$ in weak sense \cite[Proposition~5.11.1]{Kloeden1992Numerical}, 
    that is, for fixed $\mu \in \mathcal{P}([0, T) \times \R^d)$, $u \in \mathbb{U}_{\mathrm{ad}}^{\mu}$, $x = (t, z) \in \Pi_h \times \R^d$,
    \begin{equation}\label{eq_wkerror}
        \E{\varphi \circ \Phi(x, \mu, u)} = \E{\varphi(X_{t+h}^{\mu, u}) \,\big\vert\, X_t^{\mu, u} = x} + \mO(h^2) \sptext{as} h \to 0
    \end{equation}
    for any smooth test function $\varphi: [0, T) \times \R^d \to \R$, where $\circ$ denotes function composition. 
    For higher-order approximations, see, for example, the literature on numerical methods for mean-field SDEs, such as \cite{Sun2020Explicit,Sun2018Explicit,Sun2017Ito,Sun2018New}.
    \end{remark}

\begin{remark}[Clarification of the iteration index]
    The subscript $i$ in \eqref{eq_Xnm} and \eqref{eq_mun} denotes the policy iteration index in \eqref{eq_dpp_v} and \eqref{eq_dpp_u}, not a time-step index.
    Accordingly, for each fixed $i$, the samples $\{X_i^m\}_{m=1}^M$ are treated as time-space samples on $[0, T) \times \R^d$ for approximating $\mu^*$, rather than as samples at a single time level $t \in [0, T)$.
\end{remark}

\subsection{Time-discretized policy iteration}

To derive a time-discrete approximation of the dynamic programming residual $\mathcal{R}_h$ in \eqref{eq_defMh}, 
we apply the left-rectangle rule to the integral term and use the random mapping $\Phi$ to approximate the conditional distribution of $X_{(t + h)-}^{\mu, u}$ given $X_{t}^{\mu,u} = x$, yielding
\begin{equation}\label{eq_Rh_M}
    \mathcal{R}_h(x; u, v, \mu) = \E{\mathcal{M}(x; \mu, u, v)} + \mO(h^2) \quad \text{as } h \to 0+,
\end{equation}
for any $x \in \Pi_h \times \R^d$, where $\mathcal{M}$ is the martingale increment function defined by
\begin{equation}\label{eq_Mx}
    \mathcal{M}(x; \mu, u, v) := \big( h f(x, \mu, u) + v \circ \Phi(x, \mu, u) - v(x) \big) I(t \in \Pi_h), \;\; x = (t, z) \in \bar{\Pi}_h \times \R^d. 
\end{equation}
The expectation in \eqref{eq_Rh_M} is taken with respect to the Gaussian random variable $\zeta$ in \eqref{eq_randmap}. 
Under sufficient regularity conditions, the $\mO(h^2)$ term in \eqref{eq_Rh_M} arises from two sources: the left-rectangle quadrature error and the weak local truncation error of $\Phi$ described in \eqref{eq_wkerror}.
The factor $I(t \in \Pi_h)$ in \eqref{eq_Mx} is included for technical convenience, because the PE/PI steps in \eqref{eq_dpp_v} and \eqref{eq_dpp_u} are activated only for $t \in [0, T)$.

To obtain a time-discrete analogue of the PE and PI steps in \eqref{eq_dpp_v} and \eqref{eq_dpp_u}, we proceed as follows:
(i) replace $\mathcal{R}_h$ in \cref{eq_dpp_v,eq_dpp_u} by the approximation \eqref{eq_Rh_M} and neglect the $\mO(h^2)$ term;
(ii) enforce the resulting conditions at the sampled states $\{X_i^m\}_{m \in \mathbb{A}_i}$ updated at iteration $i$ according to \eqref{eq_Xnm};
(iii) replace the equilibrium measure $\mu^*$ by its empirical counterpart $\mu_i$, and approximate the pair $(u, v)$ by the neural networks $(u_{\alpha}, v_{\theta})$.
This yields the following time-discrete PE and PI steps: at iteration $i$, for each $m \in \mathbb{A}_i$,
\begin{alignat}{4}
    &\text{train } v_{\theta} &\text{ s.t. }& \E{\mathcal{M}(X_i^m; \mu_i, u_{\alpha}, v_{\theta}) \mid X_i^m} = 0, \label{eq_pe} \\
    &\text{train } u_{\alpha} &\text{ s.t. }& \min_{\alpha} \E{\mathcal{M}(X_i^m; \mu_i, u_{\alpha}, v_{\theta}) \mid X_i^m}, \label{eq_pi}
\end{alignat}
where $\mathcal{M}$ is defined in \eqref{eq_Mx}.

\begin{remark}[Connections to fictitious play]
The policy iteration in \eqref{eq_pe} and \eqref{eq_pi} is closely connected to fictitious play for mean-field games.
Specifically, at iteration $i$, the control network $u_{\alpha}$ is updated using the current empirical measure $\mu_i$ as the population distribution; in the terminology of fictitious play, this corresponds to the representative agent's best response to the population distribution inherited from the previous iteration.
This idea has been widely adopted in related works; see, for example, \cite{Cardaliaguet2017Learning,Han2020Deep,Hu2021Deep,Han2024Learning,Liang2024Actor}.
Our method differs from existing approaches in that it updates the measure via a one-step Markovian transition in \eqref{eq_Xnm}, rather than by simulating the entire state trajectory over $[0, T)$ at each iteration.
This feature enables parallelization across the time dimension and thus improves the computational efficiency of the algorithm.
\end{remark}

\subsection{Policy iteration in weak form}

The PE/PI steps above involve conditional expectations, which are difficult to evaluate accurately because, for $m \in \mathbb{A}_i$, \eqref{eq_Xnm} provides only one sample of the next state $X_{i+1}^m$ for each current state $X_i^m$.
To avoid conditional expectations, we reformulate the PE and PI steps in \eqref{eq_pe} and \eqref{eq_pi} into the following integrated form: at iteration $i$, for each $m \in \mathbb{A}_i$, 
\begin{alignat}{4}
    &\text{train } v_{\theta} &\text{ by } &\min_{\theta} \sup_{\rho \in \mathbb{T}} \big\vert \E{\rho(X_i^m) \mathcal{M}(X_i^m; \mu_i, u_{\alpha}, v_{\theta})} \big\vert^2, \label{eq_Wmartv} \\
    &\text{train } u_{\alpha} &\text{ by } &\min_{\alpha} \E{\mathcal{M}(X_i^m; \mu_i, u_{\alpha}, v_{\theta})}, \label{eq_Wmartu}
\end{alignat}
where $\mathbb{T}$ is a suitable class of test functions, for example,
\begin{equation}\label{eq_rho}
    \mathbb{T} := \bbr{\rho: [0, T) \times \R^d \to [-1, 1]^r \;\big\vert\; \rho \text{ is smooth}}.
\end{equation}
We next explain how \eqref{eq_Wmartv} and \eqref{eq_Wmartu} relate to the original PE and PI steps in \eqref{eq_pe} and \eqref{eq_pi}.

The PE step in \eqref{eq_Wmartv} can be viewed as a Galerkin-type reformulation of \eqref{eq_pe}. 
Indeed, by the tower property of conditional expectation \cite[Theorem~8.14]{Achim2020Probability},
\begin{equation*}
    \E{\rho(X_i^m)\,\mathcal{M}(X_i^m; \mu_i, u_{\alpha}, v_{\theta})} = \mathbb{E}\big[\rho(X_i^m)\,\E{\mathcal{M}(X_i^m; \mu_i, u_{\alpha}, v_{\theta}) \mid X_i^m}\big],\quad \forall \rho \in \mathbb{T}.
\end{equation*}
Therefore, \eqref{eq_Wmartv} in essence seeks to train $v_{\theta}$ such that 
\begin{equation}\label{eq2_EM}
    \mathbb{E}\big[\rho(X_i^m)\,\E{\mathcal{M}(X_i^m; \mu_i, u_{\alpha}, v_{\theta}) \mid X_i^m}\big] = 0, \quad \forall \rho \in \mathbb{T}.
\end{equation}
The above equality is precisely a weighted Galerkin formulation of the equality in \eqref{eq_pe}: the inner conditional expectation is the residual expected to vanish, while the outer expectation is a weighted inner product, testing this residual against functions in $\mathbb{T}$ with the law of $X_i^m$ as the weighting measure.

The minimization in \eqref{eq_Wmartu} is an averaged version of \eqref{eq_pi} with respect to the law of $X_i^m$. 
This reformulation avoids pointwise optimization over all the sampled states in $\{X_i^m\}_{m=1}^M$ while preserving pointwise optimality for each state in the support of the law of $X_i^m$.
Actually, by the law of total expectation, the objective in \eqref{eq_Wmartu} can be expressed as
\begin{equation*}
    \E{\mathcal{M}(X_i^m; \mu_i, u_{\alpha}, v_{\theta})} = \mathbb{E}\big[ \E{\mathcal{M}(X_i^m; \mu_i, u_{\alpha}, v_{\theta}) \mid X_i^m} \big].
\end{equation*}
Assume that the control network $u_{\alpha}$ is sufficiently expressive. 
If $u_{\alpha}$ attains the minimum of the left-hand side of \eqref{eq_Wmartu}, then the random variable $u_{\alpha}(X_i^m)$ minimizes the inner conditional expectation on the right-hand side for almost every realization of $X_i^m$, that is, for almost every state $x$ in the support of the law of $X_i^m$.
Hence, the averaged formulation in \eqref{eq_Wmartu} preserves the pointwise optimality required by the original PI step \eqref{eq_pi}. 
Similar ideas have been used in \cite[subsection~3.1]{Al2022Extensions} for solving HJB equations, and a rigorous justification is given in \cite[Lemma~3.2]{cai2024socmartnet}.

\begin{remark}[Trick of vector-valued test function]
    As indicated in \eqref{eq_rho}, the test function $\rho$ is $\R^r$-valued. 
    Therefore, in \eqref{eq_Wmartv}, the expectation is also $\R^r$-valued, and $|\cdot|$ 
    is interpreted as the Euclidean norm in $\R^r$.
    Compared with the scalar case $r=1$, using vector-valued test functions, for example with $r \geq 600$, can substantially improve the stability of the adversarial training for solving \eqref{eq_Wmartv}, as observed in \cite{cai2024socmartnet,Cai2026Deep} for solving PDEs. 
\end{remark}

\begin{remark}[Inspiration from weak adversarial networks]
    The reformulation of the PE step from \eqref{eq_pe} to \eqref{eq_Wmartv} is inspired by the weak adversarial network method for solving PDEs, originally proposed in \cite{Zang2020Weak} and subsequently developed in \cite{cai2024martingale,Cai2026Deep} for high-dimensional 
    parabolic PDEs, and stochastic control problems.
\end{remark}

\subsection{Network structure}\label{sec_net}

Now we specify the structures of the neural networks $u_{\alpha}$, $v_{\theta}$, and $\rho_{\eta}$ used in the PE and PI steps in \eqref{eq_Wmartv} and \eqref{eq_Wmartu}.
By \eqref{eq1_defUad}, the output of $u_{\alpha}$ should be restricted in the control space $U$.
Particularly, when $U = [a, b] := \prod_{j=1}^m [a_j, b_j]$ with $a_j, b_j$ the $j$-th elements of $a, b \in \R^m$, the structure of $u_{\alpha}$ can be
\begin{equation}\label{eq_defualp}
    u_{\alpha, j}(x) = \min\{\max\{a_j, \psi_{\alpha, j}(x)\}, b_j\},  \quad j = 1, 2, \cdots, m,
\end{equation}
where $\psi_{\alpha}: [0, T) \times \R^d \to \R^m$ is a neural network with parameter $\alpha$.
The network $v_{\theta}$ should fulfill the terminal condition in \eqref{eq_Vmu}, which can be achieved by the following structure: at iteration $i$,
\begin{equation}\label{eq_defvnn}
    v_{\theta}(x) := \left\{\begin{aligned}
    \phi_{\theta}(x), \quad &x \in [0, T) \times \R^d,\\ 
    g(x, \mu_i), \quad &x \in \{T\} \times \R^d,
\end{aligned}\right.
\end{equation}
where $\phi_{\theta}: [0, T) \times \R^d \to \R$ is a neural network parameterized by $\theta$, and $\mu_i$ is the empirical measure defined in \eqref{eq_mun}.

The adversarial network $\rho_{\eta}$ plays the role of test functions in classical Galerkin methods for solving PDEs.
As shown in the numerical experiments in \cite{cai2024socmartnet,Cai2026Deep}, $\rho_{\eta}$ is not necessarily very deep, but instead, it can be a shallow network with enough output dimensionality.
Following the multiscale neural network ideas in \cite{Liu2020Multi},
we consider a typical $\rho_{\eta}$ given by
\begin{equation}\label{eq_defrhonet}
    \rho_{\eta}(x) = \sin \br{\Lambda \br{W x + b}} \in \R^r, 
\end{equation}
where $\eta := (W, b) \in \R^{r \times (1+d)} \times \R^r$ is the trainable parameter; $\Lambda(\cdot)$ is a scale layer defined by 
\begin{equation}\label{eq_defLamb}
    \Lambda(y_1, y_2, \cdots, y_r) = \br{c_1 y_1, c_2 y_2, \cdots, c_r y_r}^{\top} \in \R^r, \quad c_j := 1 + (j - 1) c
\end{equation}
for $y_j \in \R$, where $c > 0$ is a non-trainable hyperparameter; 
$\sin(\cdot)$ is the activation function applied to the outputs of $\Lambda$ in an element-wise manner.

\subsection{Learning algorithm}

The weak-form policy iteration steps \eqref{eq_Wmartv} and \eqref{eq_Wmartu} depend only on unconditional expectations.
These expectations can be naturally approximated by empirical averages over the samples $\{X_i^m\}_{m \in \mathbb{A}_i}$, yielding mini-batch training schemes.
Specifically, we parameterize the test function $\rho$ by a neural network $\rho_{\eta}$ with parameter $\eta$.
Define the mini-batch estimator at iteration $i$ by
\begin{equation}\label{eq_LA}
    L(u_{\alpha}, v_{\theta}, \rho_{\eta}, \mu_i; \mathbb{A}_i)
    := \frac{1}{\abs{\mathbb{A}_i}} \sum_{m \in \mathbb{A}_i} \rho_{\eta}(X_i^m)\, \mathcal{M}(X_i^m; \mu_i, u_{\alpha}, v_{\theta}),
\end{equation}
where $\abs{\mathbb{A}_i}$ denotes the cardinality of $\mathbb{A}_i$; $\mathcal{M}$ is defined in \eqref{eq_Mx}.
Then we have the mini-batch estimators for the objectives in \eqref{eq_Wmartv} and \eqref{eq_Wmartu}:
\begin{align}
    \big\vert \E{\rho_{\eta}(X_i^m)\,\mathcal{M}(X_i^m; \mu_i, u_{\alpha}, v_{\theta})} \big\vert^2 &\approx L(u_{\alpha}, v_{\theta}, \rho_{\eta}, \mu_i; \mathbb{A}_{i, 1}) L^{\top}(u_{\alpha}, v_{\theta}, \rho_{\eta}, \mu_i; \mathbb{A}_{i, 2}), \label{eq_appr_pe} \\
    \E{\mathcal{M}(X_i^m; \mu_i, u_{\alpha}, v_{\theta})} &\approx L(u_{\alpha}, v_{\theta}, 1, \mu_i; \mathbb{A}_{i}), \quad \mathbb{A}_{i} = \mathbb{A}_{i, 1} \cup \mathbb{A}_{i, 2} \label{eq_appr_pi}
\end{align}
where $\mathbb{A}_{i, 1}$ and $\mathbb{A}_{i, 2}$ form a random partition of $\mathbb{A}_{i}$ introduced in \eqref{eq_Xnm}.
The symbol ``$1$'' in \eqref{eq_appr_pi} denotes the constant function equal to one.

Using the approximations in \eqref{eq_appr_pe} and \eqref{eq_appr_pi}, the PE step \eqref{eq_Wmartv} can be implemented via gradient descent on the right-hand side of \eqref{eq_appr_pe} with respect to $\theta$ and gradient ascent with respect to $\eta$.
Similarly, the PI step \eqref{eq_Wmartu} is implemented via gradient descent on the right-hand side of \eqref{eq_appr_pi} with respect to $\alpha$.
The samples $\{X_i^m\}_{m \in \mathbb{A}_i}$ and the empirical measure $\mu_i$ are updated according to \eqref{eq_Xnm} and \eqref{eq_mun}, respectively, using the current network $u_{\alpha}$ obtained from the PI step.
Iterating the above three steps yields our overall algorithm for solving the regenerative MFG problem.
\Cref{alg_mini} summarizes the algorithmic details.

\begin{remark}[Unbiasedness of the mini-batch approximation]
    Assuming that the sample size $M$ is sufficiently large so that the randomness of $\mu_i$ can be neglected, the disjointness of $\mathbb{A}_{i,1}$ and $\mathbb{A}_{i,2}$ implies that $L(u_{\alpha}, v_{\theta}, \rho_{\eta}, \mu_i; \mathbb{A}_{i,1})$ and $L(u_{\alpha}, v_{\theta}, \rho_{\eta}, \mu_i; \mathbb{A}_{i,2})$ are approximately independent.
    Consequently, the approximation in \eqref{eq_appr_pe} is nearly unbiased. 
    Unbiasedness of the mini-batch approximation is favorable for the convergence of stochastic gradient algorithms; see \cite{Hu2025Bias}.
\end{remark}

\begin{algorithm}[t]
    \caption{Policy iteration for the mean-field game~\eqref{eq_equilibrium}}\label{alg_mini}
    \begin{algorithmic}[1]
        \Require 
        $v_{\theta}$/$u_{\alpha}$/$\rho_{\eta}$: neural network parameterized by $\theta$/$\alpha$/$\eta$;
        $\delta_1$/$\delta_{2}$/$\delta_{3}$: learning rate for the network $v_{\theta}$/$u_{\alpha}$/$\rho_{\eta}$;
        $I$: maximum number of iterations;
        $J$/$K$: number of $(\theta, \alpha)$/$\eta$ updates per iteration;
        $M$: full batch size;
        $\hat{\mu}$: initial guess of the equilibrium measure $\mu^*$;
        
        \State Initialize the networks $u_{\alpha}$, $v_{\theta}$ and $\rho_{\eta}$

        \State Generate i.i.d. samples $\{X_0^{m}\}_{m=1}^M$ from $\hat{\mu}$  
        
        \For{$i = 0, 1, \cdots, I$}
        
        \State Generate $\mu_i$ as the empirical measure of $\{X_i^m\}_{m=1}^M$ by \eqref{eq_mun}
        
        \State Sample disjoint index sets $\mathbb{A}_{i, 1}$ and $\mathbb{A}_{i, 2}$ from $\{1, 2, \cdots, M\}$
        
        \For{$j = 0, 1, \cdots ,J-1$} \Comment{$L$ is given in \eqref{eq_LA}}
        \State $\theta \leftarrow \theta - \delta_1 \nabla_{\theta} \bbr{L(u_{\alpha}, v_{\theta}, \rho_{\eta}, \mu_i; \mathbb{A}_{i, 1}) L^{\top}(u_{\alpha}, v_{\theta}, \rho_{\eta}, \mu_i; \mathbb{A}_{i, 2})}$
        \State $\alpha \leftarrow \alpha - \delta_2 \nabla_{\alpha} L(u_{\alpha}, v_{\theta}, 1, \mu_i; \mathbb{A}_{i, 1} \cup \mathbb{A}_{i, 2})$
        \EndFor
        \For{$k = 0, 1, \cdots, K-1$}
        \State $\eta \leftarrow \eta + \delta_{3} \nabla_{\eta} \bbr{L^{\top}(u_{\alpha}, v_{\theta}, \rho_{\eta}, \mu_i; \mathbb{A}_{i, 1}) L(u_{\alpha}, v_{\theta}, \rho_{\eta}, \mu_i; \mathbb{A}_{i, 2})}$
        \EndFor
        
        \For{$m = 1, 2, \cdots, M$} \Comment{$\Phi$ is given in \eqref{eq_randmap}}
        \State $X_{i+1}^m = \Phi(X_i^m, \mu_i, u_{\alpha})$ if $m \in \mathbb{A}_{i, 1} \cup \mathbb{A}_{i, 2}$; otherwise, $X_{i+1}^m = X_i^m$
        \EndFor
        \EndFor
        \Ensure $u_{\alpha}$, $v_{\theta}$ and $\mu_I$
    \end{algorithmic}
\end{algorithm}

We close this section with several remarks on the practical advantages of the proposed method.
\begin{itemize}
    \item \textbf{Avoiding higher-order derivatives:} 
    Related neural approaches, such as \cite{Assouli2024Deep,Ruthotto2020Machine,Lin2021Alternating}, are built around the coupled HJB--FP structure, or closely related value--density formulations. Such approaches often involve second-order spatial derivatives and, in some cases, Hessian-related computations. 
    Automatic differentiation for such PDE terms becomes expensive when the state dimension $d$ is large \cite{He2023Learning,hu2024sdgd,Hu2024Hutchinson,shi2024stochastic,Cai2026Deep}.
    Our formulation avoids the explicit computation of these derivatives by performing policy iteration in weak form.
    
    \item \textbf{Avoiding full path simulation:}
    In many trajectory-based approaches to MFGs and related stochastic games or McKean--Vlasov problems, e.g., \cite{Han2020Deep,Hu2021Deep,Han2024Learning,Liang2024Actor},
    the distributional component is approximated by simulating multiple trajectories of $X^{\mu, u}$.
    This requires generating sample paths sequentially in time, which can be computationally expensive when the time step size $h$ is small.
    By contrast, the regenerative formulation uses the random mapping $\Phi$ to update samples through one-step transitions, thereby avoiding full path simulation at each iteration.

    \item \textbf{Avoiding pointwise Hamiltonian minimization:} 
    The HJB--FP system contains an infimum over the control space $U$ in the HJB equation.
    This infimum may not admit a closed-form solution and must therefore be approximated numerically, leading to pointwise optimization over the state space.
    Such computations become particularly expensive when $U$ is high-dimensional.
    In the proposed method, these pointwise optimization problems are replaced by an averaged optimization problem for training the control network $u_{\alpha}$.
    
    \item \textbf{Shared computation results:}
    The sample ensemble $\{X_i^m\}_{m=1}^M$ supplies the spatio-temporal samples used in the PE and PI steps.
    These samples reflect the state dynamics and concentrate on the most relevant regions of the state space.
    In addition, the ensemble update via \eqref{eq_Xnm} and the PE/PI steps in \eqref{eq_Wmartv} and \eqref{eq_Wmartu} all use the same transitions, namely $\{X_i^m\}_{m \in \mathbb{A}_i} \to \{X_{i+1}^m\}_{m \in \mathbb{A}_i}$.
    Thus, the PE/PI, measure update, and network training steps share a common set of sample transitions, reducing the overall computational cost.
\end{itemize}
Taken together, these features make the proposed method well suited to high-dimensional MFG, as illustrated by the numerical examples in the next section.

\section{Numerical examples}\label{sec_num}

The numerical examples in this section are designed to demonstrate the performance of \Cref{alg_mini} for a series of MFG problems.
Throughout, the abbreviations ``RE'', ``RC'', ``SD'', ``RT'', ``MEM'', ``vs'' stand for ``relative error'', ``relative cost'', ``standard deviation'', ``runtime'', ``gpu memory usage'', ``versus'', respectively. 

The parameter settings of \Cref{alg_mini} are as follows.
The neural networks $u_{\alpha}$ and $v_{\theta}$ are constructed as in \eqref{eq_defualp} and \eqref{eq_defvnn}, respectively, with fully connected inner networks $\psi_{\alpha}$ and $\phi_{\theta}$ of depth $H = 6$ and width $W$, where $W = 104$ for $d = 1$ and $W = 1008$ for $d = 1000$.
Both networks use ReLU activation functions.
The adversarial network $\rho_{\eta}$ is chosen as in \eqref{eq_defrhonet}, with output dimension $r = 1200$ and scaling parameter $c = 10$ in \eqref{eq_defLamb}.
At iteration $i$, the learning rates are set to $\delta_1 = \delta_2 = \delta_0 \times 10^{-3} \times 0.01^{i/I}$ for $u_{\alpha}$ and $v_{\theta}$, and $\delta_0 \times 10^{-2} \times 0.01^{i/I}$ for $\rho_{\eta}$, with $\delta_0 := 3d^{-0.5}$.
The maximum number of iterations is set to $I = 9,\!000$, 
and the inner loop numbers are set to $J = 2K = 2$.
The full batch size is set to $M = 1024 \times 10^3$, and the mini-batch sizes are set to $\abs{\mathbb{A}_{i,1}} = \abs{\mathbb{A}_{i,2}} = M/20$ for all $i$.
The initial guess $\hat{\mu}$ of the equilibrium measure $\mu^*$ is specified in each example below.
The time step size in \eqref{eq_Pih} is set to $h = T/100$, and stochastic gradient updates are performed by RMSProp.
All experiments are implemented in Python 3.12 using PyTorch 2.6.0, with float32 precision.
The networks $u_{\alpha}$ and $v_{\theta}$ are trained using the Automatic Mixed Precision technique\footnote{\url{https://docs.pytorch.org/docs/2.7/amp.html}}.
The algorithms are accelerated using Distributed Data Parallel\footnote{\url{https://github.com/pytorch/tutorials/blob/main/intermediate_source/ddp_tutorial.rst}} on a compute node equipped with 8 NVIDIA GeForce RTX 4090 GPUs.

The initial state $X_0 = (0, Z_0)$ in \eqref{eq_state} is uniformly distributed on the diagonal of a $d$-dimensional hypercube; that is, $Z_0 \sim \mathrm{Uniform}(D)$ with 
\begin{equation}\label{eq_D}
    D := \{s \B{1}_d: s \in [-c, c]\}, \quad \B{1}_d := (1, 1, \cdots, 1)^{\top} \in \R^d. 
\end{equation}
We take $c = 3$ in \cref{sec_sr} and $c = 1$ in \cref{sec_numerical_lq}.
The performance of the method is evaluated by comparing $J(u_{\alpha}, \hat{\mu}_I)$ with the theoretical optimal cost, and by comparing $v_{\theta}$ with the exact $v$. 
The following metrics quantify this comparison:
\begin{equation*}
    \mathrm{RC} := \frac{\hat{J} - J^*}{J^*}, \quad \mathrm{RE}_1 := \frac{\sum_{z \in D_{\mathrm{re}}} |v_{\theta}(0, z) - v(0, z)|}{\sum_{z \in D_{\mathrm{re}}} |v(0, z)|}, \quad \mathrm{RE}_{\infty} := \frac{\max_{z \in D_{\mathrm{re}}} |v_{\theta}(0, z) - v(0, z)|}{\max_{z \in D_{\mathrm{re}}} |v(0, z)|}. 
\end{equation*}
Here $D_{\mathrm{rc}}$ and $D_{\mathrm{re}}$ are test point sets with $|D_{\mathrm{rc}}| = 256$ and $|D_{\mathrm{re}}| = 1,\!000$, randomly sampled from $D$ and fixed during training.

The cost $\hat{J}$ is the empirical version of $J(u_{\alpha}, \hat{\mu}_i)$ given by \eqref{eq_cost}, where the expectation is approximated by averaging over $256$ paths generated using Euler--Maruyama approximation of the state dynamics~\eqref{eq_state} with $(u, \mu)$ replaced by $(u_{\alpha}, \hat{\mu}_i)$ at iteration $i$, and the time integral is evaluated using left-rectangle quadrature with step size $h = T/100$.
The cost $J^*$ approximates the theoretical optimal cost $J(u^*, \mu^*)$ by averaging $v(0, z)$ over the test point set $z \in D_{\mathrm{rc}}$. 
Thus, $\mathrm{RC}$ quantifies the relative error in the cost of the learned control $u_{\alpha}$ compared to the theoretical optimum, while $\mathrm{RE}_1$ and $\mathrm{RE}_{\infty}$ quantify the accuracy of the learned value function $v_{\theta}$ relative to the exact $v$.

\subsection{Linear-quadratic MFG with explicit solutions}\label{sec_numerical_lq}

\begin{figure}[t]
    \centering
    \resizebox{1.0\textwidth}{!}{
    \begin{tabular}{@{}c@{\hspace{0mm}}ccc@{}}
        & Trajectory & RE vs Iter. & RC vs Iter.   \\
        \adjustbox{valign=m}{\rotatebox[origin=c]{90}{LQ-1}} 
        & 
        \adjustbox{valign=m}{\includegraphics[width=0.35\textwidth]{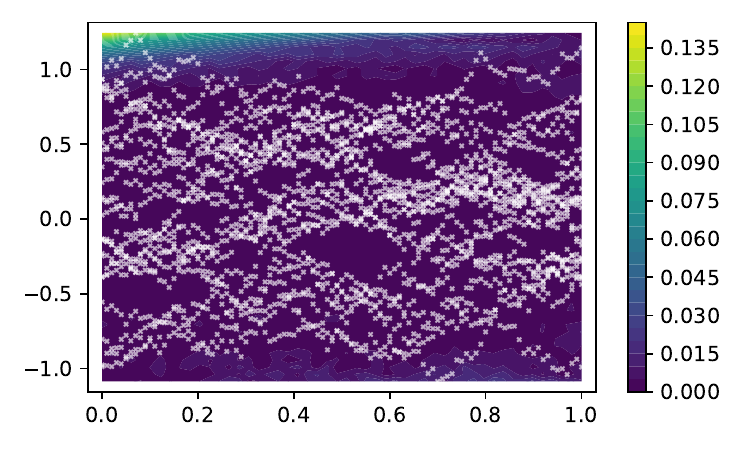}} 
        &
        \adjustbox{valign=m}{\includegraphics[width=0.27\textwidth]{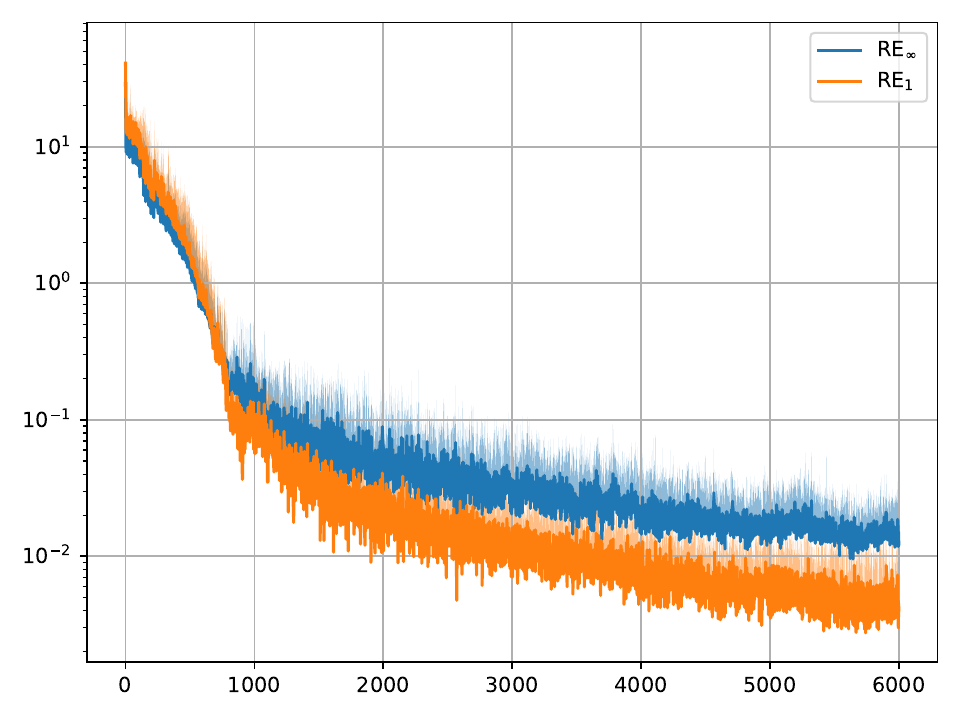}} 
        &
        \adjustbox{valign=m}{\includegraphics[width=0.27\textwidth]{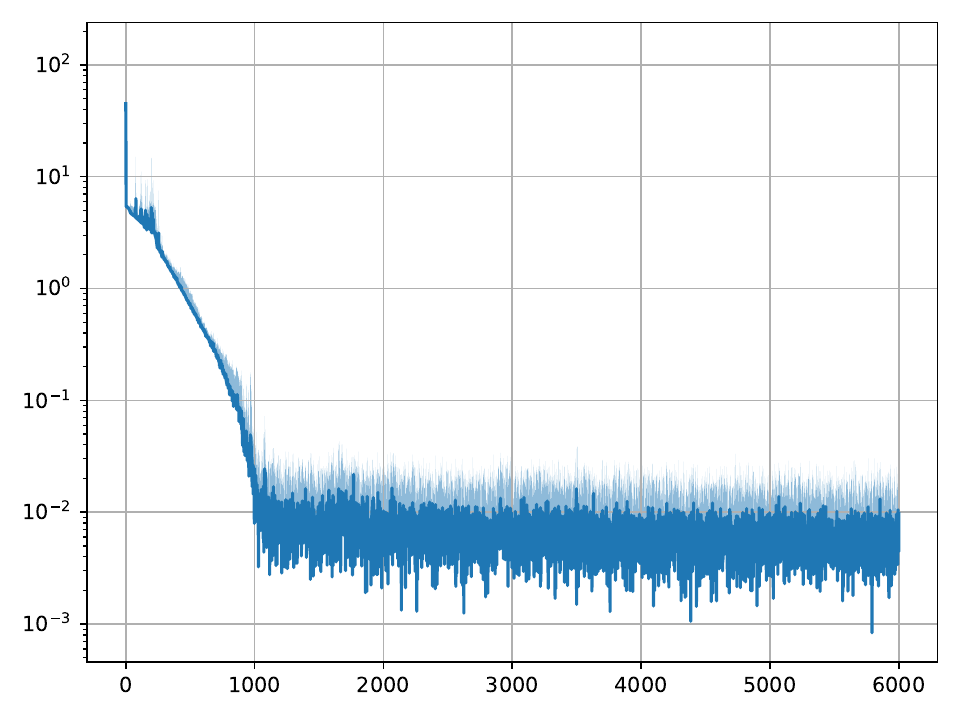}} 
        \\
        \adjustbox{valign=m}{\rotatebox[origin=c]{90}{LQ-2}} 
        & 
        \adjustbox{valign=m}{\includegraphics[width=0.35\textwidth]{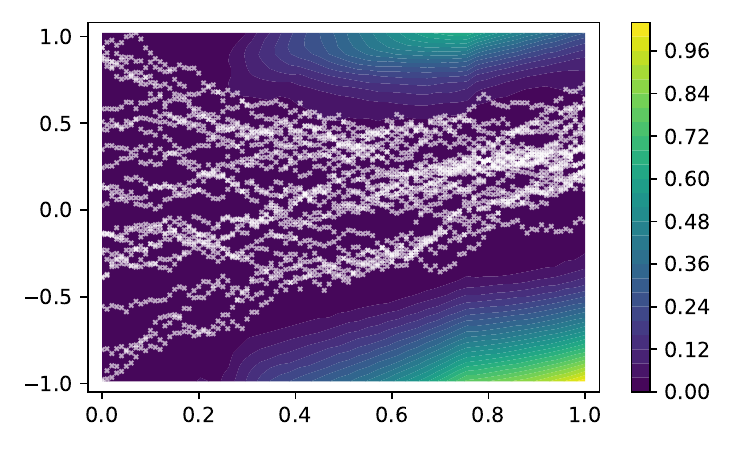}} 
        &
        \adjustbox{valign=m}{\includegraphics[width=0.27\textwidth]{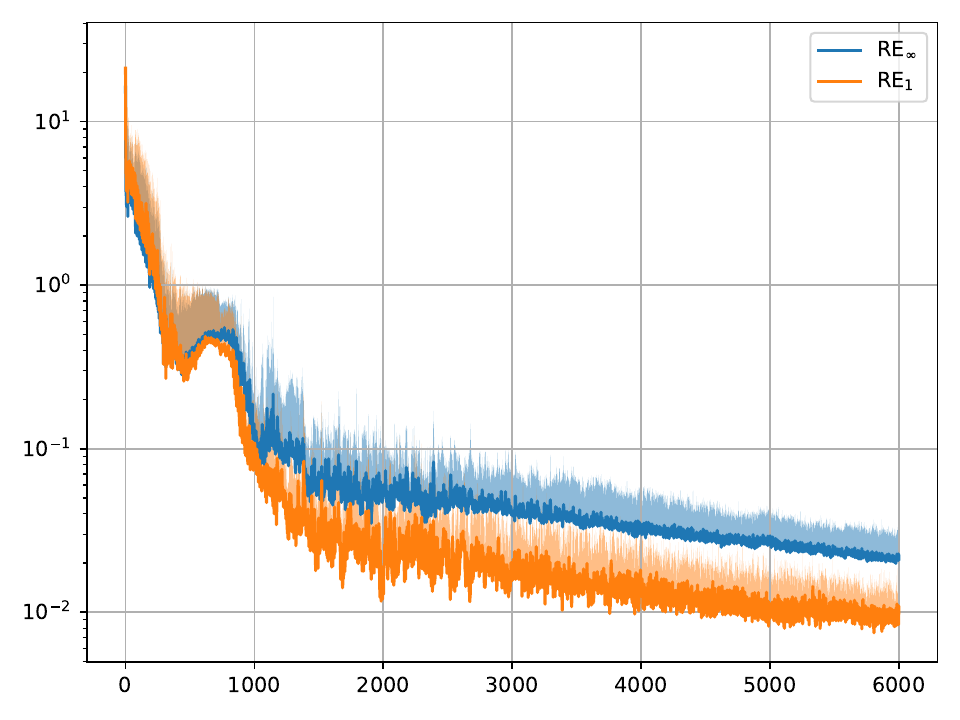}} 
        &
        \adjustbox{valign=m}{\includegraphics[width=0.27\textwidth]{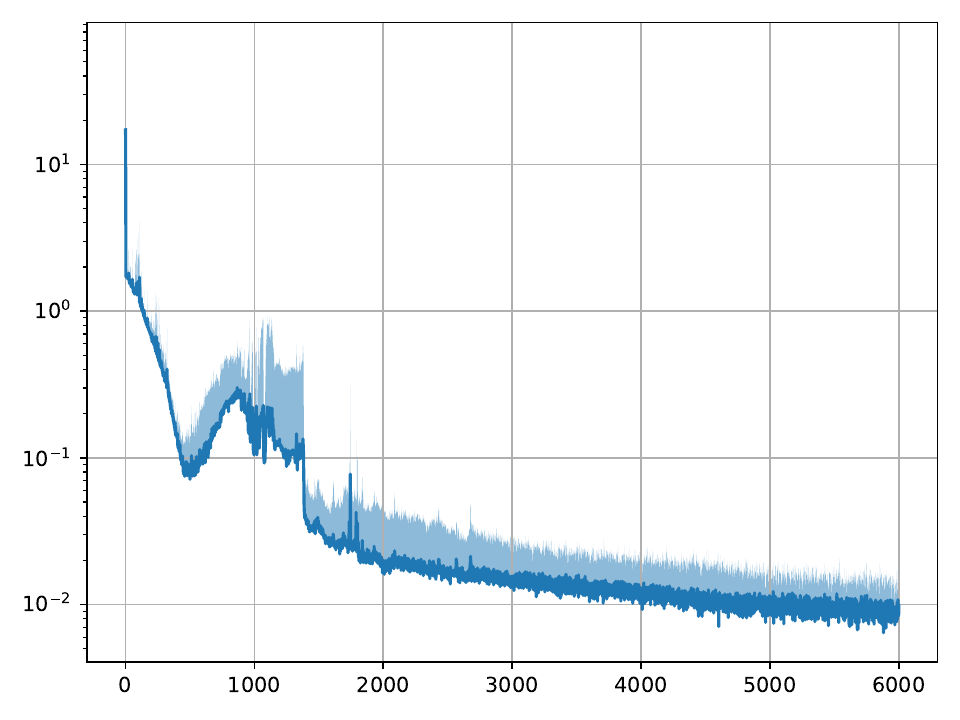}} 
        \\
        \adjustbox{valign=m}{\rotatebox[origin=c]{90}{LQ-3}} 
        & 
        \adjustbox{valign=m}{\includegraphics[width=0.35\textwidth]{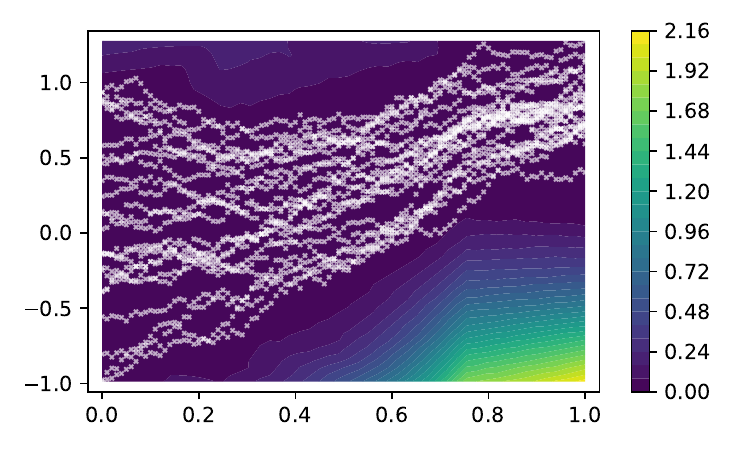}} 
        &
        \adjustbox{valign=m}{\includegraphics[width=0.27\textwidth]{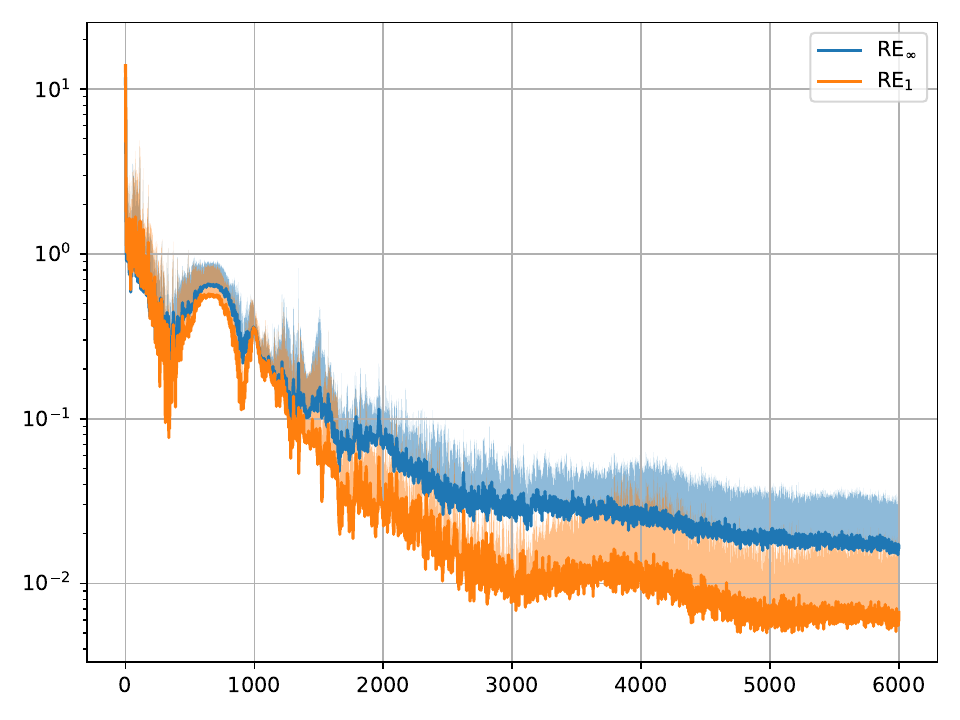}} 
        &
        \adjustbox{valign=m}{\includegraphics[width=0.27\textwidth]{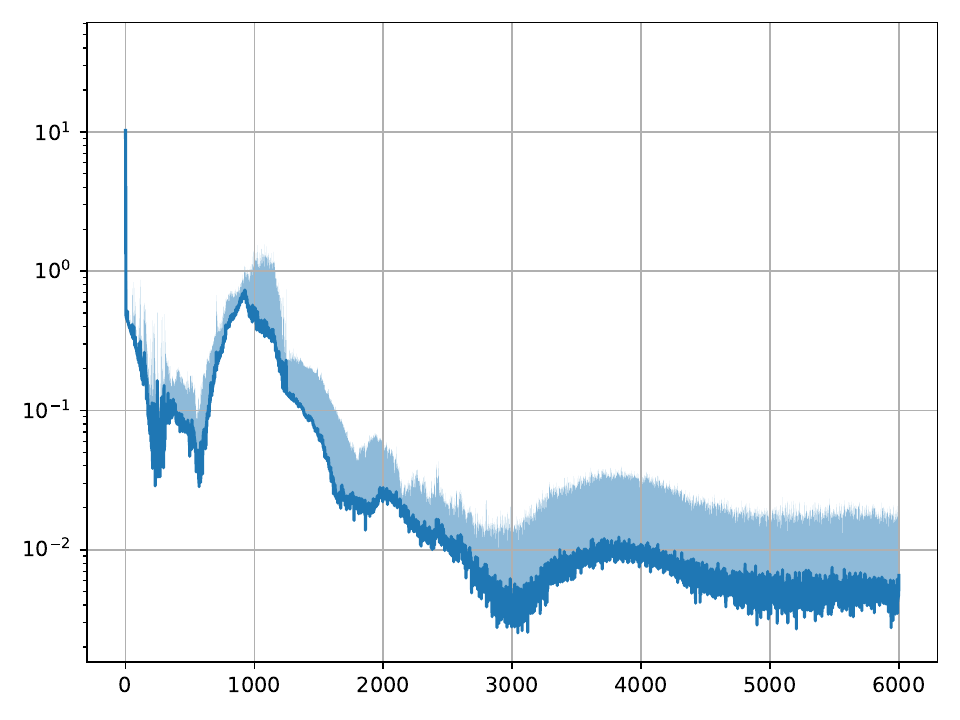}} 
    \end{tabular}
    }
    \caption{
    Numerical results of \Cref{alg_mini} for LQ-1, -2, and -3 in \cref{sec_numerical_lq} with $d = 1$. 
    In the first column, scattered points show simulated paths of $X_t^{\mu_I, u_{\alpha}}$ from \eqref{eq_Xnm}; the background color shows the pointwise absolute error of the value network $v_{\theta}$ at the final iteration; the x- and y-axes denote time $t$ and state $z$, respectively.
    In the second and third columns, the shaded regions indicate the mean plus $2 \times \text{SD}$ over 5 independent runs.
    \Cref{tab_ReSdRtHjb} summarizes the final quantitative results.
}\label{fig_LQd1}
\end{figure}

\begin{figure}[t]
    \centering
    \resizebox{1.0\textwidth}{!}{
    \begin{tabular}{@{}c@{\hspace{1mm}}cccc@{}}
        & \textbf{$s \mapsto v(0, s\B{1}_d)$} & \textbf{RE vs Iter.} & \textbf{RC vs Iter.} & \textbf{Loss vs Iter.} \\
        
        \adjustbox{valign=m}{\rotatebox[origin=c]{90}{LQ-1}} & 
        \adjustbox{valign=m}{\includegraphics[width=0.36\textwidth]{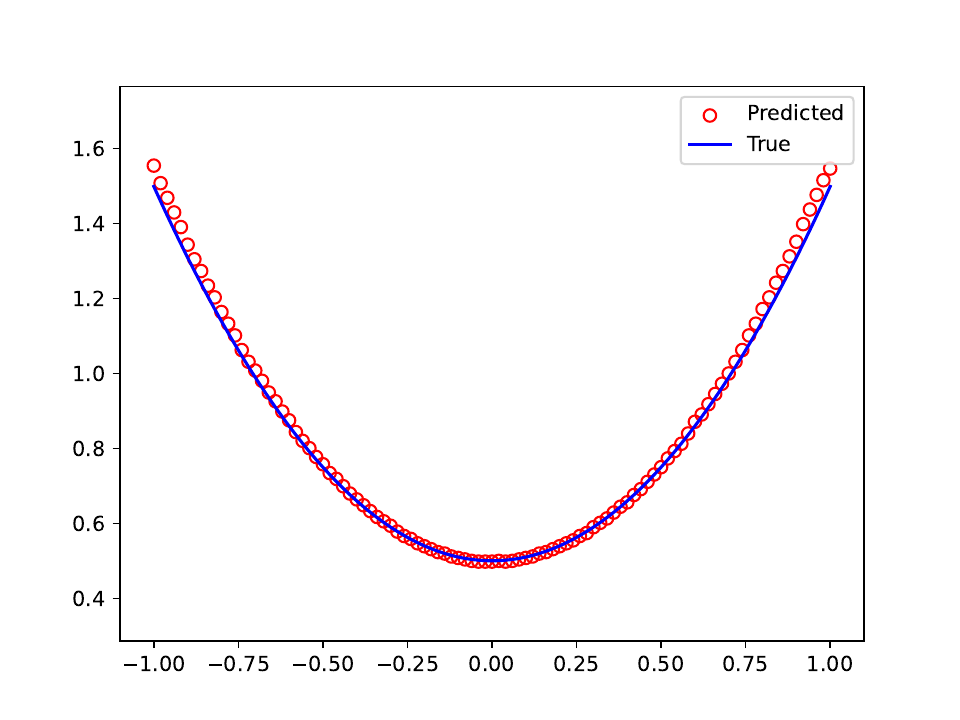}} &
        \adjustbox{valign=m}{\includegraphics[width=0.31\textwidth]{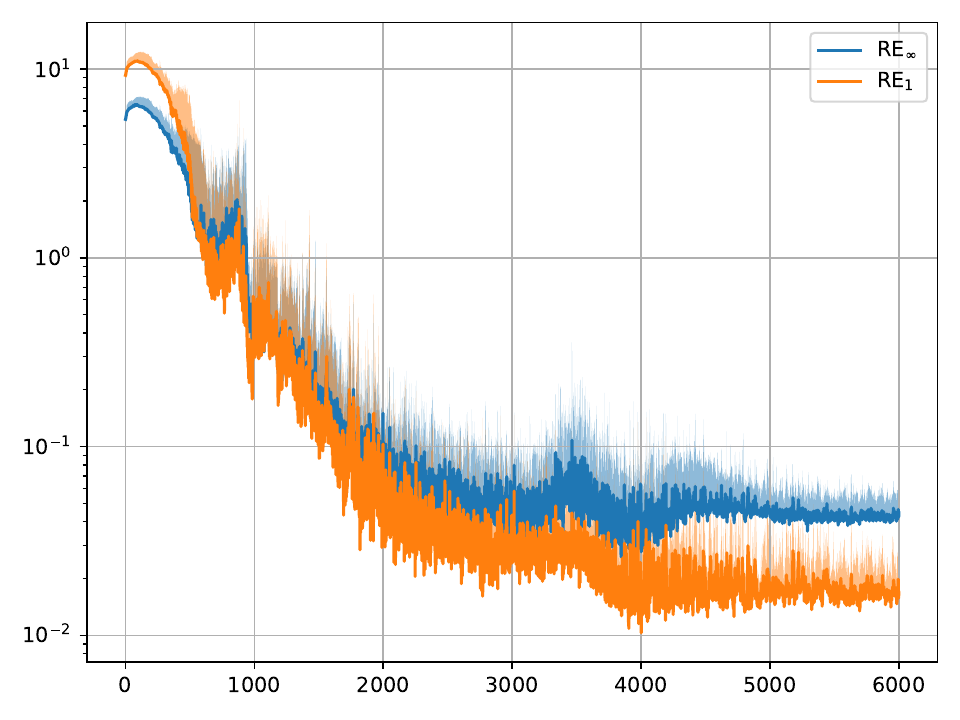}} &
        \adjustbox{valign=m}{\includegraphics[width=0.3\textwidth]{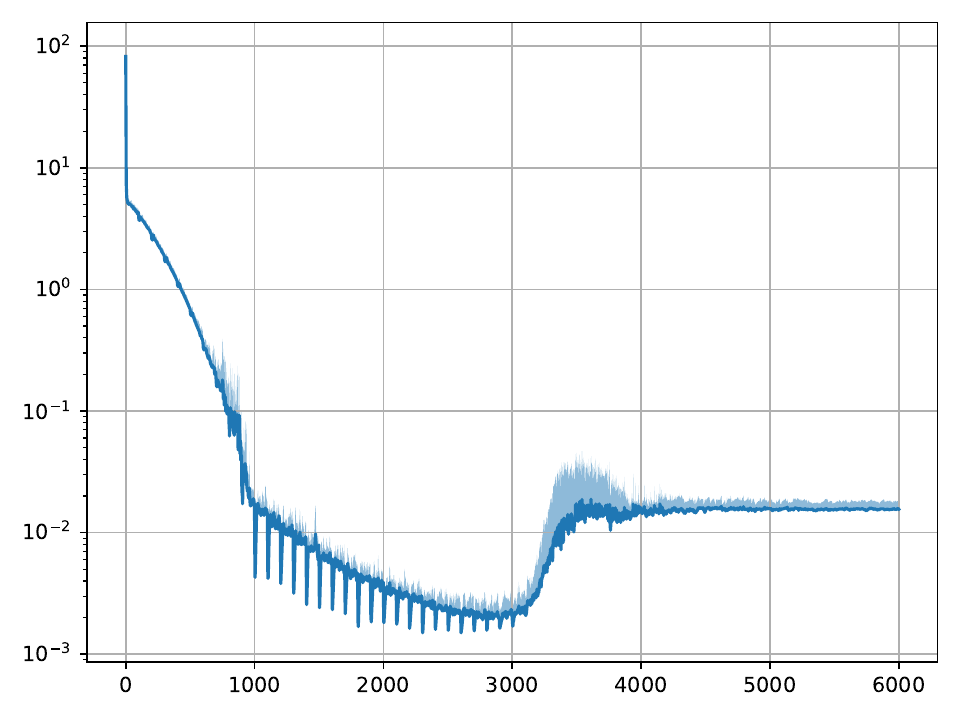}} &
        \adjustbox{valign=m}{\includegraphics[width=0.3\textwidth]{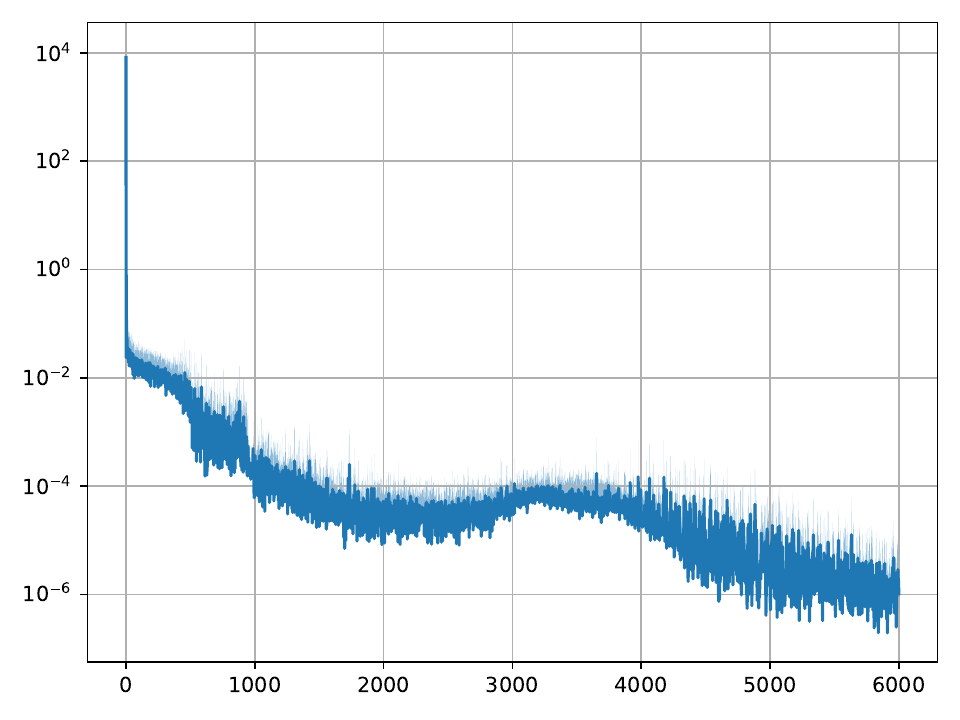}} \\
        
        \adjustbox{valign=m}{\rotatebox[origin=c]{90}{LQ-2}} & 
        \adjustbox{valign=m}{\includegraphics[width=0.36\textwidth]{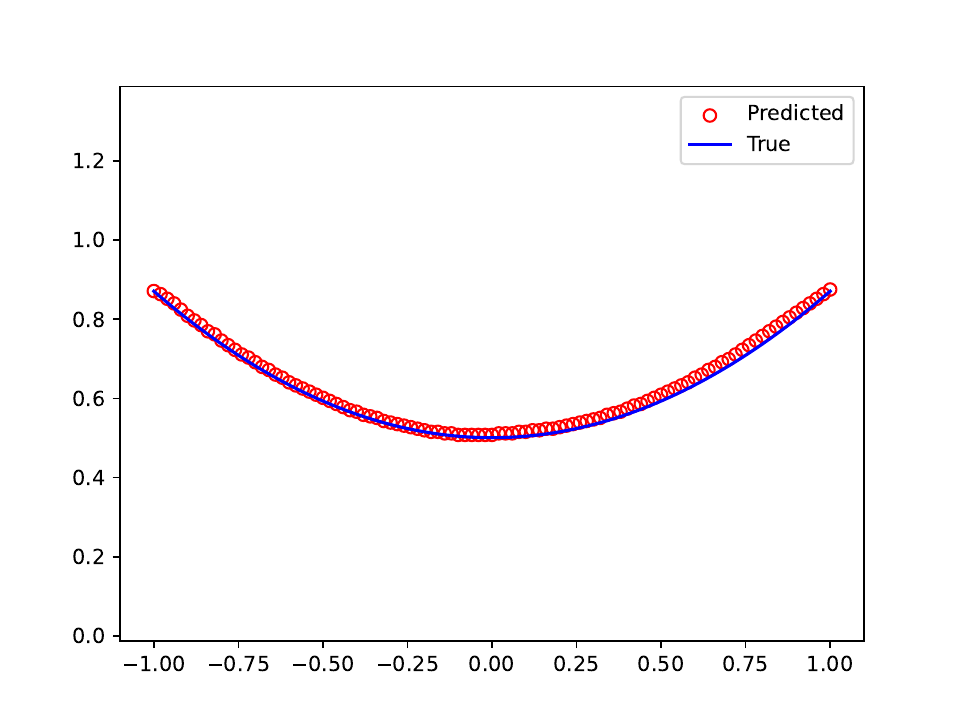}} &
        \adjustbox{valign=m}{\includegraphics[width=0.31\textwidth]{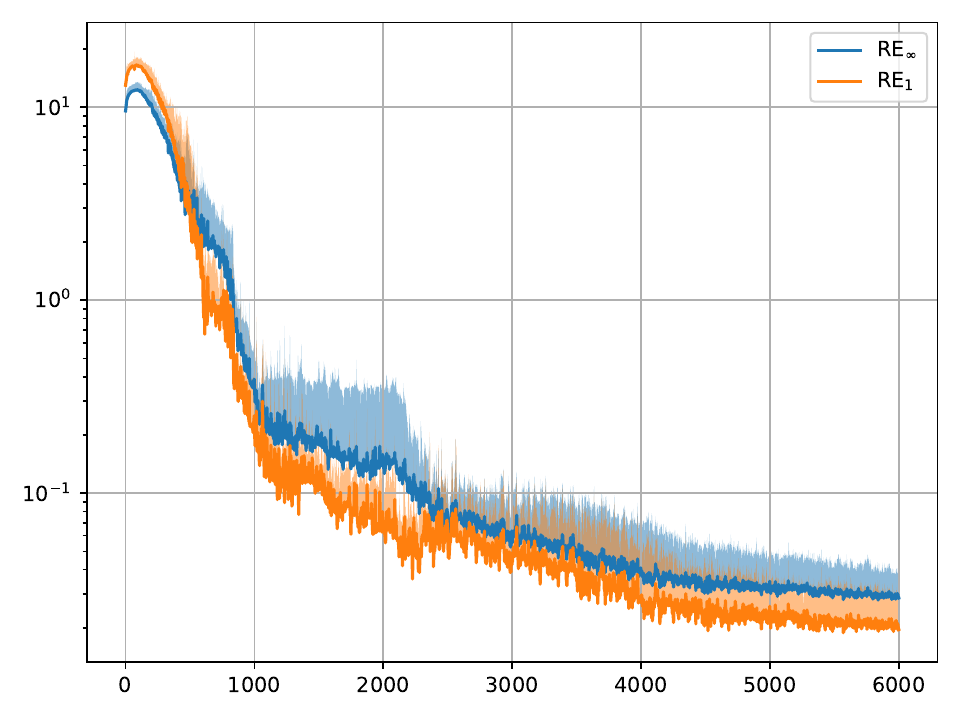}} &
        \adjustbox{valign=m}{\includegraphics[width=0.3\textwidth]{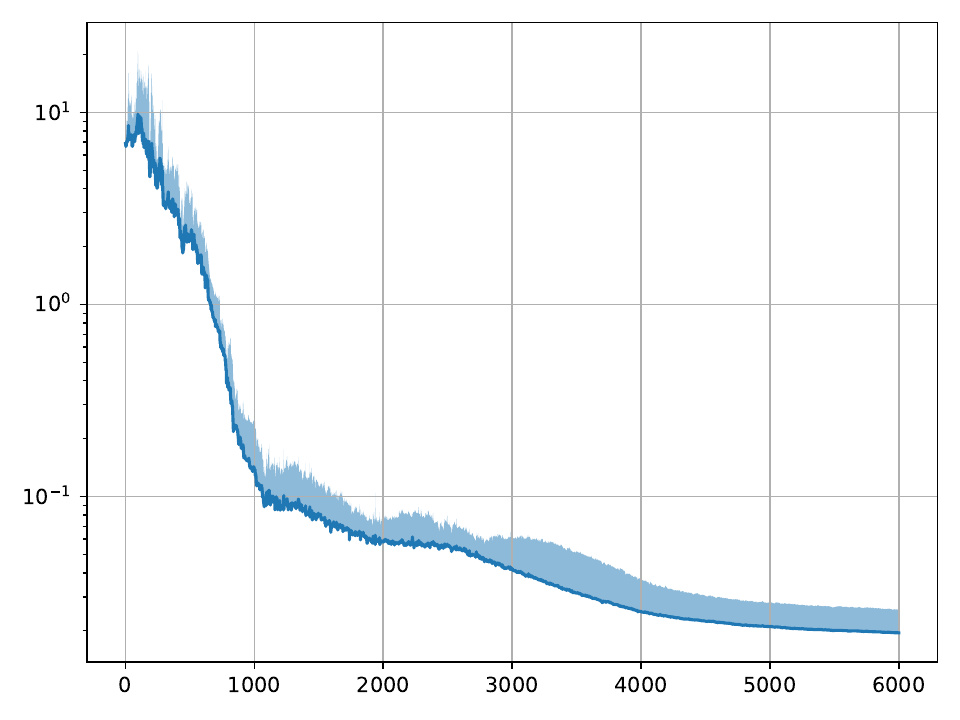}} &
        \adjustbox{valign=m}{\includegraphics[width=0.3\textwidth]{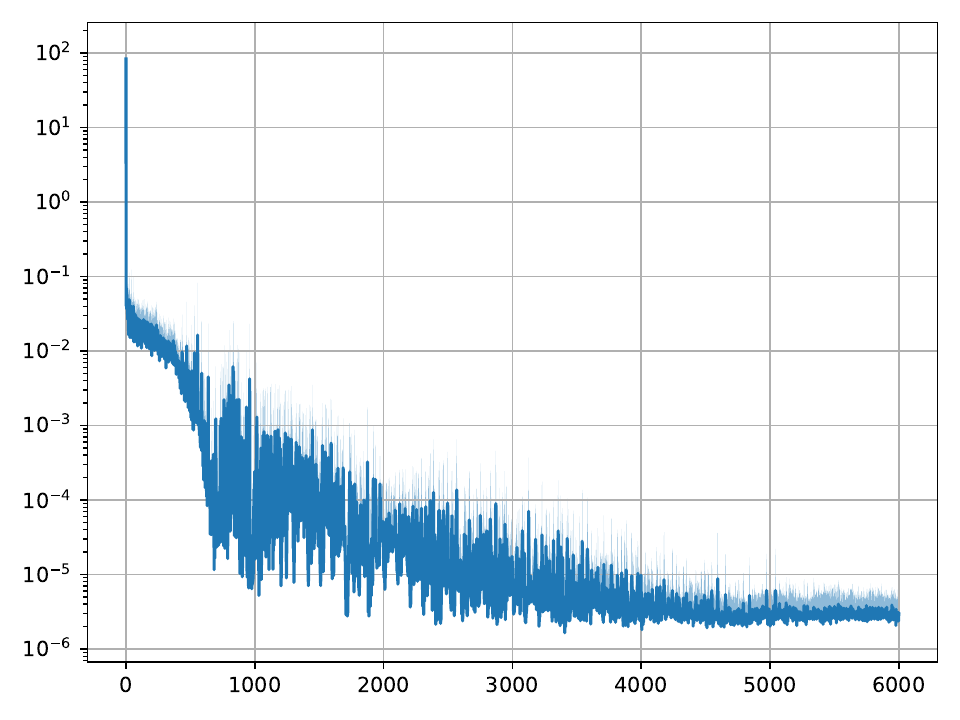}} \\
        
        \adjustbox{valign=m}{\rotatebox[origin=c]{90}{LQ-3}} & 
        \adjustbox{valign=m}{\includegraphics[width=0.36\textwidth]{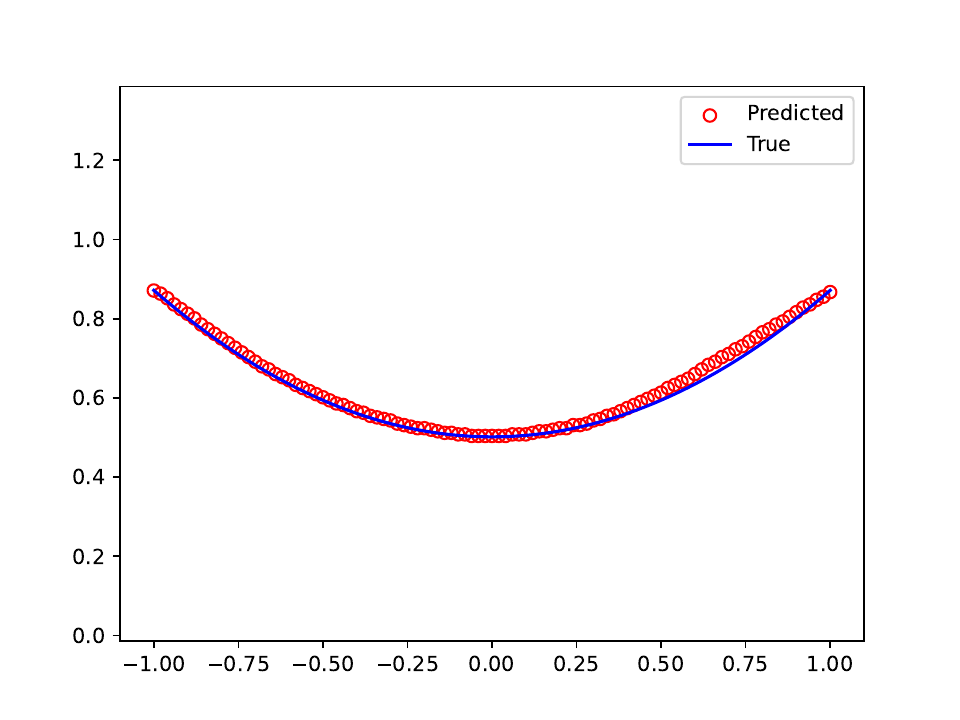}} &
        \adjustbox{valign=m}{\includegraphics[width=0.31\textwidth]{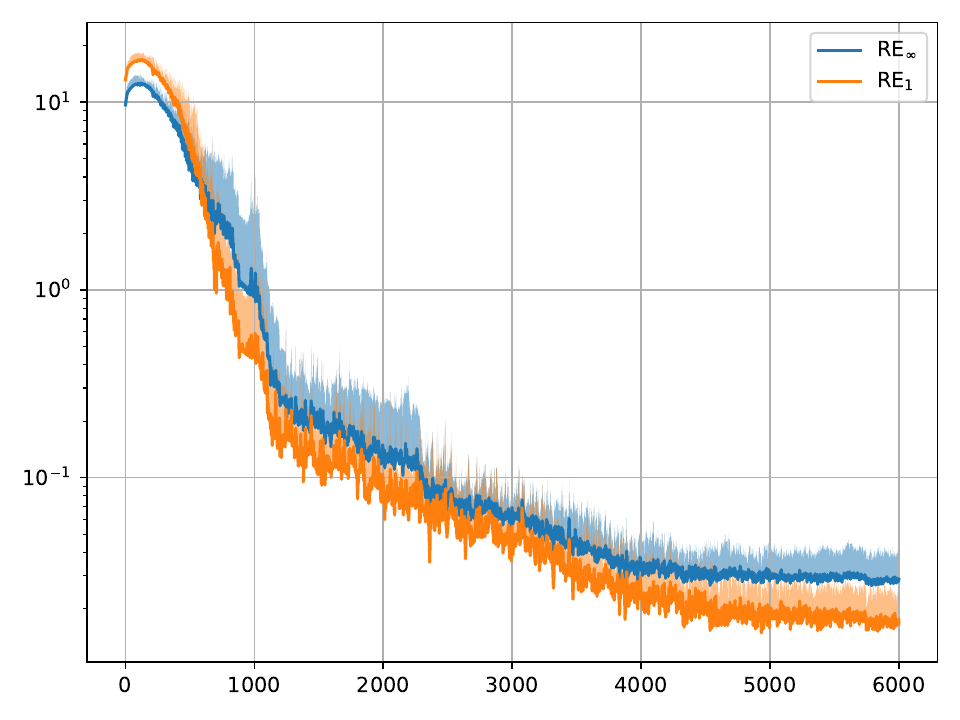}} &
        \adjustbox{valign=m}{\includegraphics[width=0.3\textwidth]{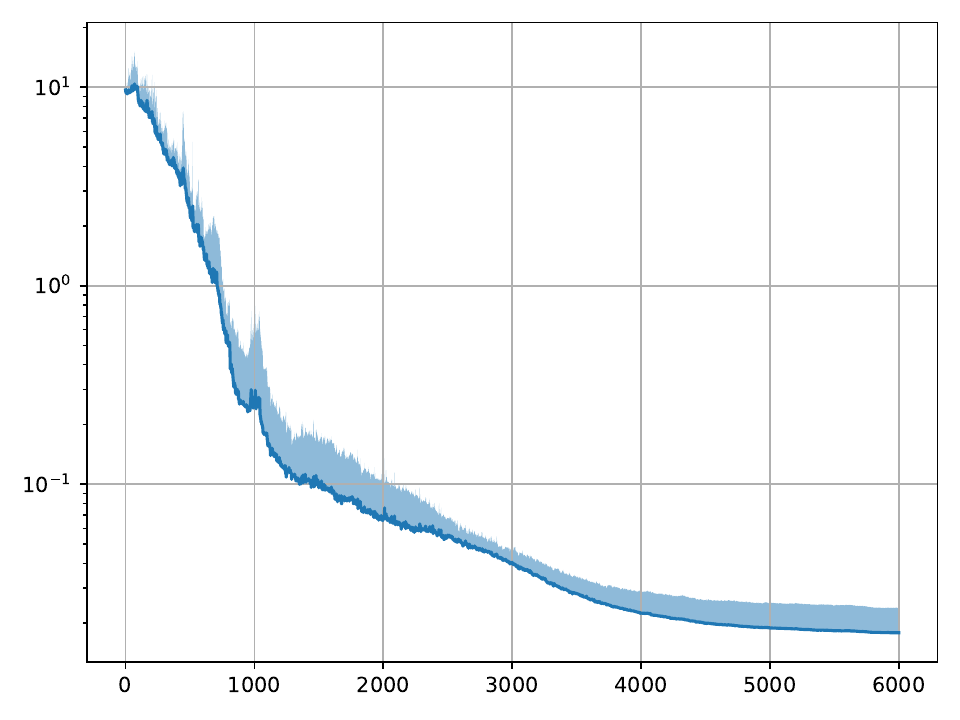}} &
        \adjustbox{valign=m}{\includegraphics[width=0.3\textwidth]{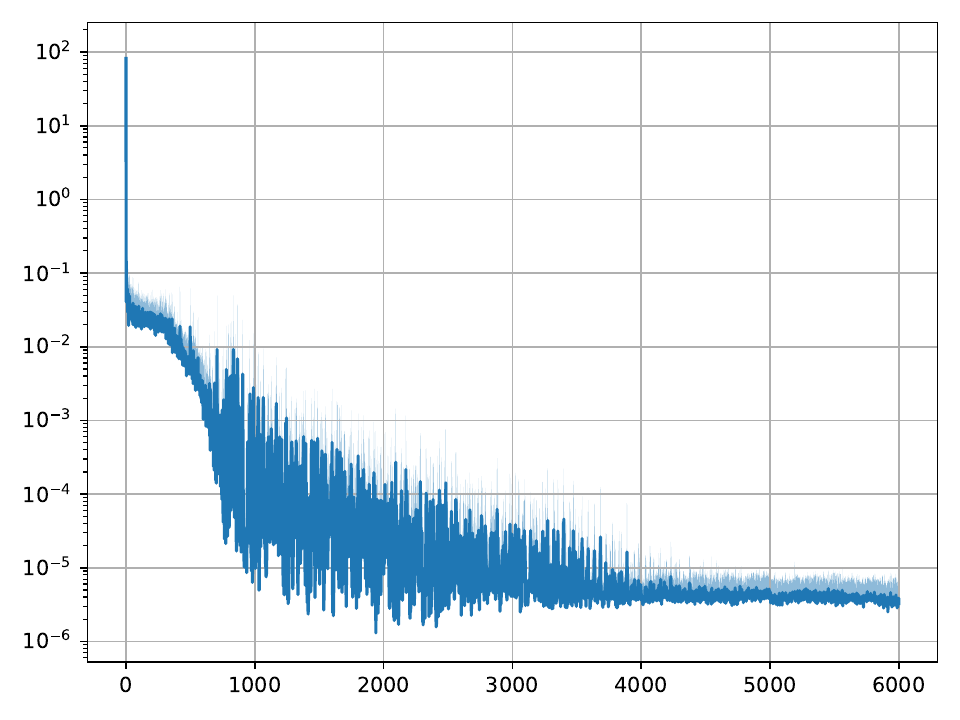}}
    \end{tabular}
    }
    \caption{
    Numerical results of \Cref{alg_mini} for LQ-1, -2, and -3 in \cref{sec_numerical_lq} with $d = 1,000$. 
    The first column shows the exact and predicted values of $s \mapsto v(0, s\B{1}_d)$. 
    In the second to fourth columns, the shaded regions show the mean plus $2 \times \text{SD}$ over five independent runs.
    Loss in the last-column heading denotes the absolute value of the right-hand side of \eqref{eq_appr_pe}.
    \Cref{tab_ReSdRtHjb} summarizes the final quantitative results.
    }\label{fig_LQd1e3}
\end{figure}

\begin{figure}[t]
    \includegraphics[width=0.33\textwidth]{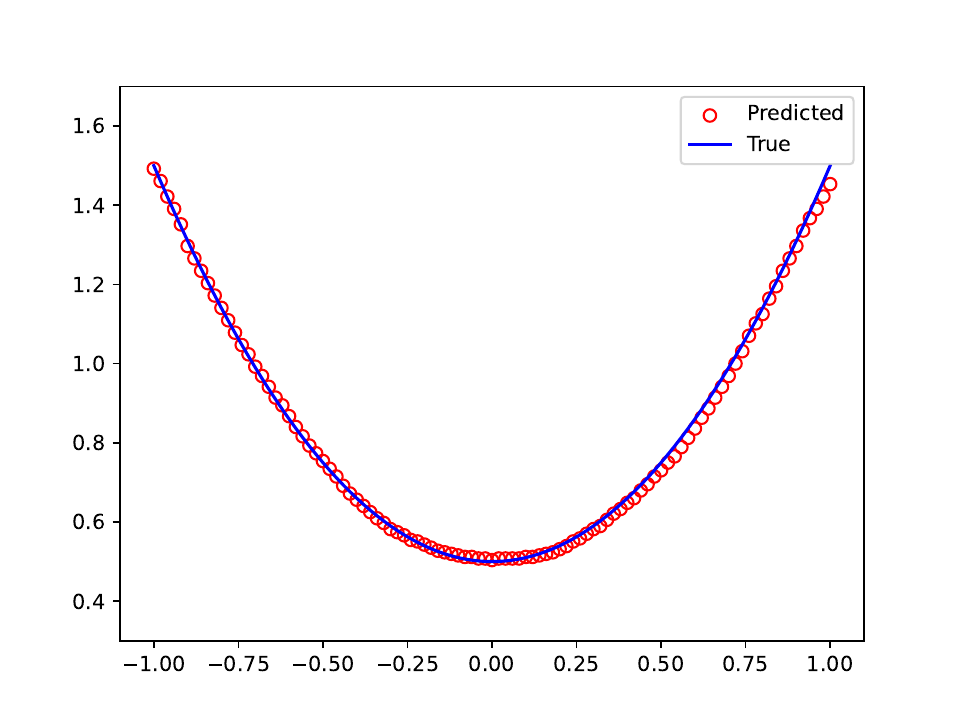}
    \includegraphics[width=0.3\textwidth]{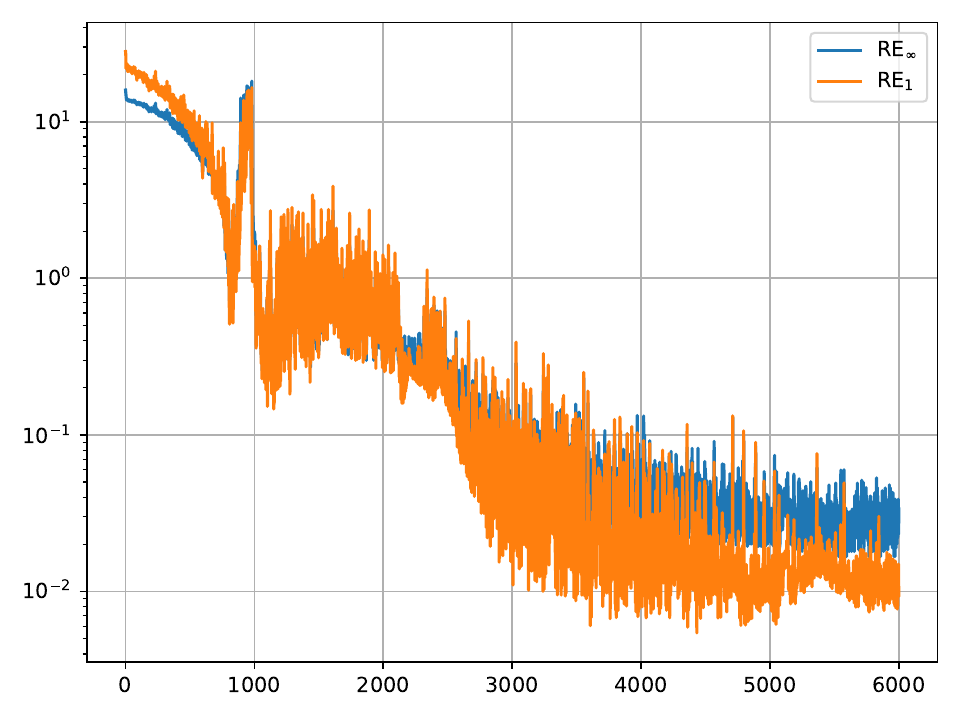}
    \includegraphics[width=0.3\textwidth]{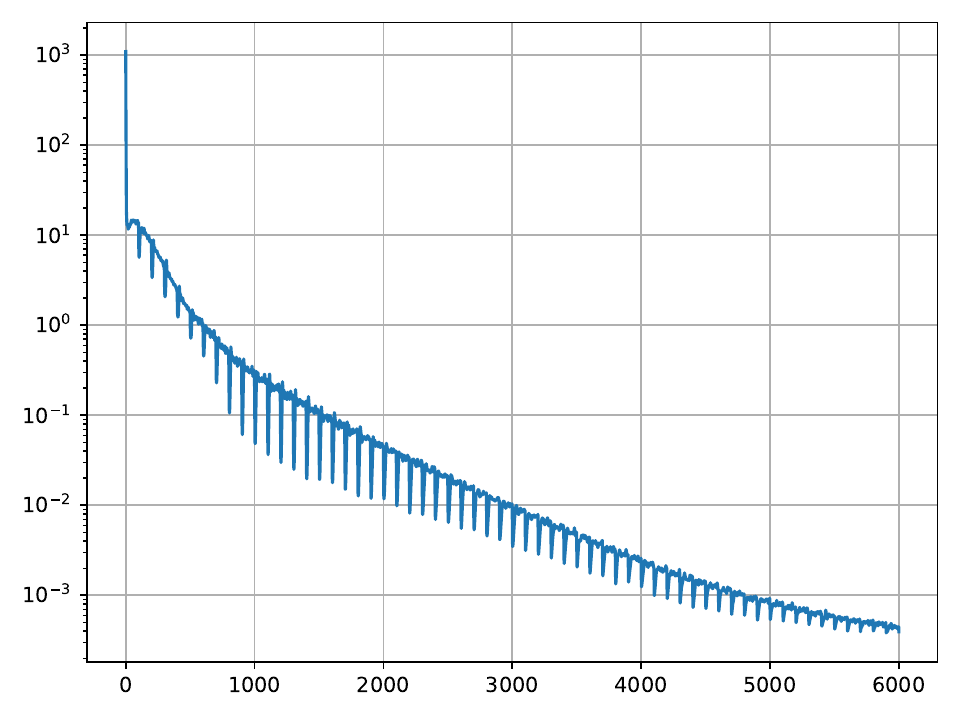}
    \caption{
    Numerical results of \Cref{alg_mini} for LQ-1 in \Cref{sec_numerical_lq} with $d = 10{,}000$.
    The left subplot compares the exact and predicted values of $s \mapsto v(0, s\B{1}_d)$.
    The middle and right subplots show RE and RC versus iteration.
    The runtime is 5258 seconds, and the peak per-iteration single-GPU memory usage is 55,645 MB, multiplied by 8 to obtain the total usage across 8 GPUs.
    }\label{fig_LQd1e4}
\end{figure}

This example considers a linear-quadratic MFG that admits an explicit solution and thus serves as a benchmark for validating the proposed method.
The state dynamics and the cost functional are given by \eqref{eq_state} and \eqref{eq_cost}, respectively, with the following specifications:
\begin{equation}\label{eq_state_lq}
    b_z(x, \mu, u) = u + c_0 (\bar{m}_t - z) + c_1 (z^*(t) - \bar{m}_t), \quad \sigma_z(x, \mu, u) = c_{\sigma} I_d, \quad T = 1,
\end{equation}
\begin{equation}\label{eq_fg_lg}
    f(x, \mu, u) = \frac{1}{2} \big( c_2 \abs{u}^2 + c_3 \abs{z - \bar{m}_t}^2 + c_4 \abs{z - z^*(t)}^2 \big), \quad g(x, \mu) = \frac{c_5}{2} \abs{z - z^*(T)}^2 + \frac{1}{2},
\end{equation}
for $x = (t, z) \in [0, T] \times \R^d$ and $u \in \R^d$, where $\bar{m}_t$ is the mean of the population of agents at time $t$, defined by
\begin{equation}\label{eq_mbar}
    \bar{m}_t := \int_{[0, T] \times \R^d} \delta_t(s) z \, \mu(\!\di s \times \di z), 
\end{equation}
with $\delta_t(\cdot)$ denoting the Dirac delta function centered at $t$.
The function $z^*: [0, T] \to \R^d$ denotes a deterministic target trajectory.
The scalar constants $c_0, c_1, \cdots, c_5, c_{\sigma} \in \R$ together with the target trajectory $z^*$ are specified below.
The initial guess $\hat{\mu}$ for the equilibrium measure $\mu^*$ in \Cref{alg_mini} is chosen as the law of $(t, Z_{\mathrm{init}, t})$ for $t \in [0, T]$, where $Z_{\mathrm{init}, t} := Z_0 + 5 t \B{1}_d + c_{\sigma} B_t$.
This choice of $\hat{\mu}$ is intentionally far from the true $\mu^*$ and therefore tests the robustness of the algorithm with respect to the initial guess.

The optimal control $u^*$ and the value function $v$ for this linear-quadratic MFG are given by
\begin{equation}\label{eq_ansatz_lq}
    u^*(x) = -c_2^{-1} \big(a(t)(z - \bar{m}_t) + b(t) \big), \;\;\; v(x) = \frac{1}{2} a(t) \abs{z - \bar{m}_t}^2 + b(t)^{\top}(z - \bar{m}_t) + \gamma(t) + \frac{1}{2},
\end{equation}
for $x = (t, z) \in [0, T) \times \R^d$, where the functions $a$, $b$, $\gamma$, and $\bar{m}$ satisfy the following system of ordinary differential equations (ODEs):
\begin{align}
    &a'(t) = c_2^{-1} a(t)^2 + 2c_0 a(t) - c_3 - c_4, \quad a(T)=c_5, \label{eq_final_ode_a}\\
    &\dot{\bar{m}}_t = - c_2^{-1} b(t) + c_1\bigl(z^*(t)-\bar{m}_t\bigr), \quad \bar{m}_0 = \E{Z_0}, \label{eq_final_ode_m}\\
    &b'(t) = c_0 b(t) - c_4\bigl(\bar{m}_t - z^*(t)\bigr), \quad b(T)=c_5\bigl(\bar{m}_T - z^*(T)\bigr), \label{eq_final_ode_b}\\
    &\gamma'(t) = -\frac{1}{2} \bbr{c_2^{-1} |b(t)|^2 + c_{\sigma}^2 d\,a(t) + c_4 |\bar{m}_t - z^*(t)|^2}, \quad \gamma(T)=\frac{c_5}{2}|z^*(T)-\bar{m}_T|^2, \label{eq_final_ode_g}
\end{align}
for $t \in [0, T)$.
The detailed derivation of \eqref{eq_ansatz_lq}--\eqref{eq_final_ode_g} is provided in \cref{appendix_lq} in Appendix. 
For reference solutions of $v$ and $u^*$, the above ODEs are solved numerically using the variable-step Runge--Kutta solver \texttt{scipy.integrate.solve\_ivp} with the ``RK45'' method and relative tolerance $10^{-8}$.\footnote{\url{https://docs.scipy.org/doc/scipy-1.17.0/reference/generated/scipy.integrate.solve_ivp.html}} 

We consider the following specific parameter settings for \eqref{eq_state_lq} and \eqref{eq_fg_lg}.
\begin{itemize}
    \item \textbf{LQ-1:} $c_0 = c_1 = c_4 = 0$, $c_2 = 1$, $c_3 = c_5 = 1/d$, and $c_{\sigma} = 0.5 /\sqrt{d}$. Here, $c_4 = 0$ implies that the agents do not have a target-tracking term in the cost functional.
    \item \textbf{LQ-2:} $c_0 = 1$, $c_1 = 0$, $c_2 = c_3 = c_4 = c_5 = 1 / d$, and $c_{\sigma} = 0.5 /\sqrt{d}$. The target trajectory is the unit circle defined by
    \begin{equation}\label{eq_target_circle}
        z^*(t) = y^*(t) / \abs{y^*(t)}, \quad  y_i^*(t) := \sin(2 \pi t + i \pi / 2), \quad i = 1, 2, \cdots, d, \quad t \in [0, T]. 
    \end{equation}
    \item \textbf{LQ-3:} 
    This setting is the same as LQ-2, except that $c_1 = - 0.5$ and the target trajectory is replaced by the helix curve defined by
    \begin{equation*}
        z^*(t) = 2 t\, y^*(t) / \abs{y^*(t)}, \quad t \in [0, T] 
    \end{equation*}
    with $y^*(t)$ defined in \eqref{eq_target_circle}.
\end{itemize}
The results of \Cref{alg_mini} for LQ-1, LQ-2, and LQ-3 with $d = 1$ are presented in \Cref{fig_LQd1}, while those for $d = 1000$ are presented in \Cref{fig_LQd1e3}.
We also perform an additional experiment for LQ-1 with $d = 10{,}000$ to further demonstrate the scalability of the method; the results are shown in \cref{fig_LQd1e4}.

The results in \Cref{fig_LQd1} may be summarized as follows.
In the first column, the background color shows that the absolute error of the value network $v_{\theta}$ remains small in regions visited by the simulated paths and becomes larger in less-explored regions.
The second and third columns show that the RE and RC decrease steadily over the iterations, indicating stable convergence.
In particular, the RC in the third column reaches the order of $10^{-2}$ for all three settings, suggesting that the learned control $u_{\alpha}$ attains a cost close to the theoretical optimum.
Together with the visualization in the first column, these results indicate that the proposed method effectively explores the parts of the state space that are most relevant to the optimal control problem and learns an accurate value function in those regions.

Now we turn to the results in \Cref{fig_LQd1e3} for $d = 1000$.
For visualization, the first column compares the exact and predicted values of $v(t, z)$ at the initial time $t = 0$ with $z$ along the diagonal of the state space.
Close agreement is observed.
The second and third columns show that both the RE and the RC decrease steadily over the iterations, similar to the behavior observed in \Cref{fig_LQd1} for $d = 1$.
These results suggest that the proposed method effectively learns an accurate value function and a near-optimal control even in the high-dimensional case $d = 1000$.

The quantitative summary in \Cref{tab_ReSdRtHjb} is consistent with these observations.
Across all six cases, both $\mathrm{RE}_1$ and RC are of order $10^{-2}$, and $\mathrm{RE}_{\infty}$ stays below $4.4 \times 10^{-2}$.
When the dimension increases from $d=1$ to $d=1000$, $\mathrm{RE}_1$, $\mathrm{RE}_{\infty}$, and RC increase only moderately, by factors ranging from about $2.0$ to $4.1$, $1.4$ to $3.4$, and $1.6$ to $2.7$, respectively, while RT and MEM rise by factors of about $2$ and $7.5$.
The corresponding standard deviations remain small relative to the means, indicating good robustness over five independent runs.
Overall, these results show that \Cref{alg_mini} scales well to high-dimensional MFG while maintaining a favorable balance among accuracy, stability, and computational cost.

We close this example with the $d = 10{,}000$ results in \Cref{fig_LQd1e4}.
The left subplot shows close agreement between the predicted value function $v_{\theta}$ and the exact $v$, while the middle and right subplots show steady decay of RE and RC, consistent with the lower-dimensional cases.
With a runtime of 5258 seconds and a peak per-iteration single-GPU memory usage of 55,645 MB, the method remains computationally feasible while maintaining good accuracy at $d = 10{,}000$.

\begin{table}[t]
    \centering
    \caption{
    Numerical results of \Cref{alg_mini} for LQ-1, -2, and -3 in \cref{sec_numerical_lq}. In the $\mathrm{RE}_1$, $\mathrm{RE}_{\infty}$, and RC columns, values outside and inside parentheses are the mean and standard deviation over 5 independent runs, respectively. 
    In the MEM column, values report the peak single-GPU memory usage per iteration, in megabytes (MB), multiplied by 8 to obtain the total usage across 8 GPUs.
    Convergence histories are shown in \Cref{fig_LQd1,fig_LQd1e3}.  
    }\label{tab_ReSdRtHjb}
    \resizebox{\textwidth}{!}{
    \begin{tabular}{l l l l l l l} 
    \toprule
    Equation & $d$ & $\mathrm{RE}_1$ & $\mathrm{RE}_\infty$ & RC & MEM (MB) & RT (s) \\ [0.5ex] 
    \midrule
    LQ-1 & 1    &  4.05E-3 (2.09E-3) & 1.26E-2 (5.20E-3)  &   9.69E-3 (6.18E-3) &  535  & 294      \\
    LQ-1 & 1000 & 1.68E-2 (2.13E-3) & 4.33E-2 (5.99E-3) & 1.55E-2 (1.26E-3) & 3964 & 579 \\[1ex]
    LQ-2 & 1    & 9.92E-3 (2.41E-3)  & 2.11E-2 (3.91E-3)  &  9.84E-3 (3.24E-3)  &  535  & 308      \\
    LQ-2 & 1000 & 1.97E-2 (4.57E-3) & 2.86E-2 (5.18E-3) & 1.95E-2 (3.12E-3) & 3995  & 602 \\[1ex]
    LQ-3 & 1    & 6.67E-3 (5.62E-3)  & 1.67E-2 (7.65E-3)  &  6.54E-3 (5.23E-3)  &  535  &  317     \\
    LQ-3 & 1000 & 1.75E-2 (2.15E-3) & 2.88E-2 (5.22E-3) & 1.79E-2 (3.01E-3) & 3997  & 615 \\
    \bottomrule
    \end{tabular}
    }
\end{table}

\subsection{Interbank systemic-risk MFG}\label{sec_sr}

We consider a mean-field game with clear financial interpretation, arising from a stylized model of interbank borrowing and lending \cite{Carmona15Mean}. 
For compatibility with the numerical scheme developed in this paper, we focus on a version of the model without common noise.
The state dynamics and the cost functional are given in \eqref{eq_state} and \eqref{eq_cost}, respectively, with the drift and diffusion in \eqref{eq_bzsgmz} specified as 
\begin{equation}\label{eq_sr_state}
    b_z(x, \mu, u) = c_1\bigl(\bar{m}_t - z\bigr) + u, \quad \sigma_z(x, \mu, u) = c_2, 
\end{equation}
and the running and terminal cost functions specified as
\begin{equation}\label{eq_sr_cost}
    f(x, \mu, u) = \frac{1}{2}u^2 - c_3\,u\bigl(\bar{m}_t-z\bigr) + \frac{c_4}{2} \bigl(\bar{m}_t - z\bigr)^2, \quad g(x, \mu) = \frac{c_5}{2}\bigl(\bar{m}_T - z\bigr)^2,  \quad T = 1,
\end{equation}
for $x = (t, z) \in [0, T] \times \R$ and $u \in \R$. 
Here the functions $b_z$, $\sigma_z$, $f$, and $g$ are scalar-valued, and  $\bar{m}_t$, defined in \eqref{eq_mbar}, denotes the population mean at time $t$ under $\mu$.
The model parameters are specified below:
\begin{equation}
    c_1 = 1, \quad c_2 = 0.5, \quad c_3 = 1, \quad c_4 = 2, \quad c_5 = 1.
\end{equation}

This model has a natural economic interpretation.
Under the above specifications, the state process in \eqref{eq_state} takes the form $X_t^{\mu, u} = (t, Z_t^{\mu, u})$, where $Z_t^{\mu, u}$ denotes the log-monetary reserve of a representative bank and $\bar{m}_t$ is the population mean reserve, representing the aggregate liquidity of the banking system.
In \eqref{eq_sr_state}, the drift term $c_1(\bar{m}_t - z)$ models interbank borrowing and lending: banks with reserves below the average are pulled upward, whereas banks with reserves above the average are pushed downward. 
Thus, $c_1$ measures the strength of the stabilizing mean reversion generated by the interbank market.
The control $u$ represents the bank's net borrowing/lending rate with the central bank, while the constant volatility $c_2$ captures idiosyncratic liquidity shocks.
In the running cost \eqref{eq_sr_cost}, the term $u^2/2$ is the direct cost of using the central bank liquidity facility, while the cross term $-c_3\,u(\bar{m}_t - z)$ encourages a bank below the average reserve to borrow and a bank above the average reserve to lend. 
The parameter $c_3$ can be interpreted as an effective incentive or fee parameter set by the regulator. 
The penalties with coefficients $c_4$ and $c_5$ penalize deviations from the system average over the time interval and at the terminal time, respectively.
Following \cite{Carmona15Mean}, this formulation reflects the idea that the central bank acts as a clearing house providing liquidity support, while systemic stress is associated with broad downward movements in banks' reserve levels.

\begin{figure}[H]
    \centering
    \includegraphics[width=0.28\textwidth]{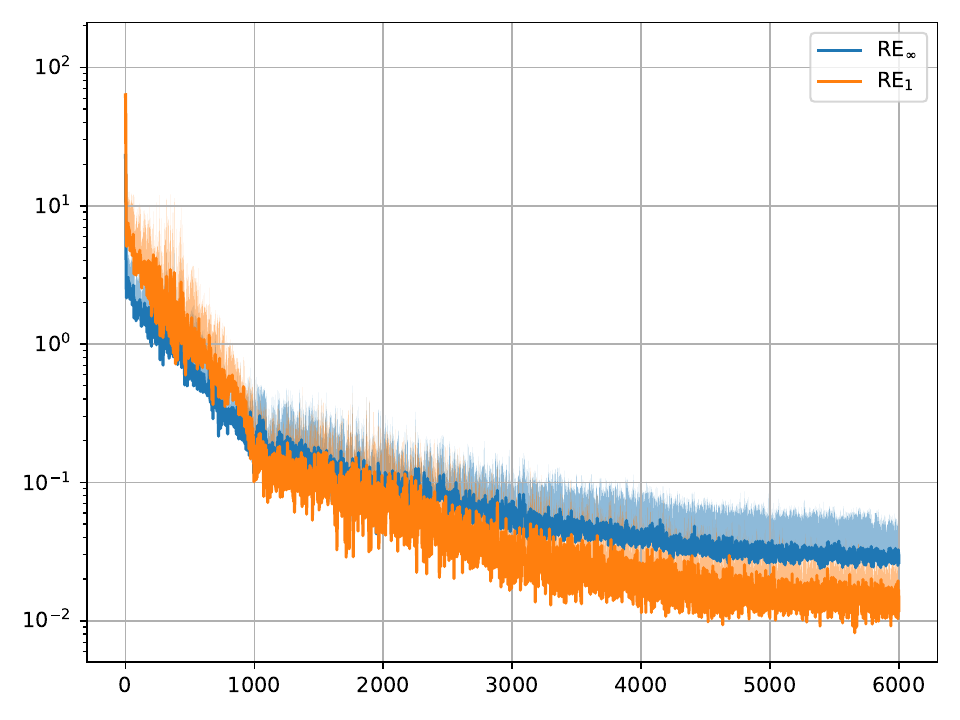}
    \includegraphics[width=0.28\textwidth]{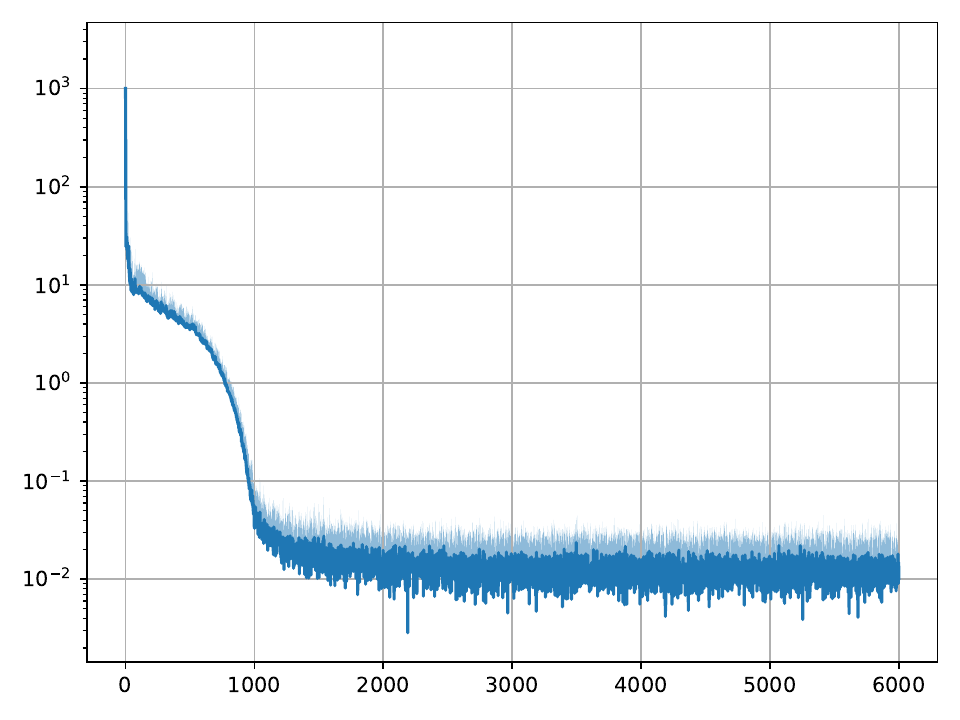}
    \includegraphics[width=0.35\textwidth]{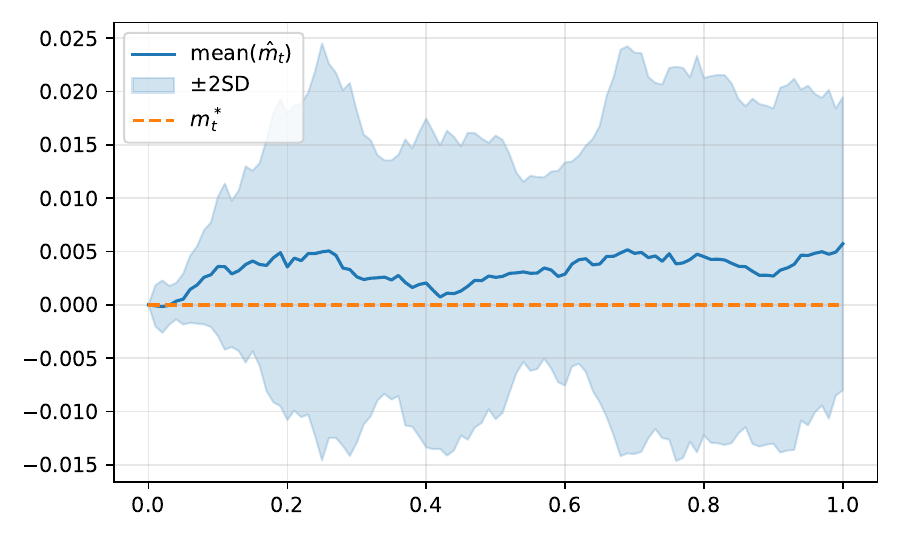}
    \\
    \includegraphics[width=0.95\textwidth]{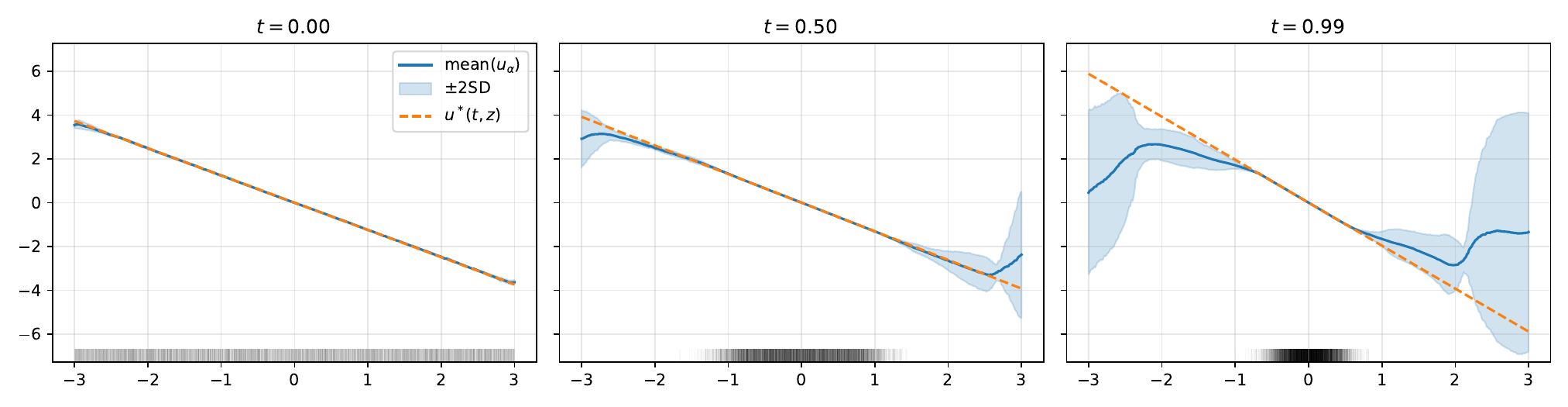}\\
    \includegraphics[width=0.4\textwidth]{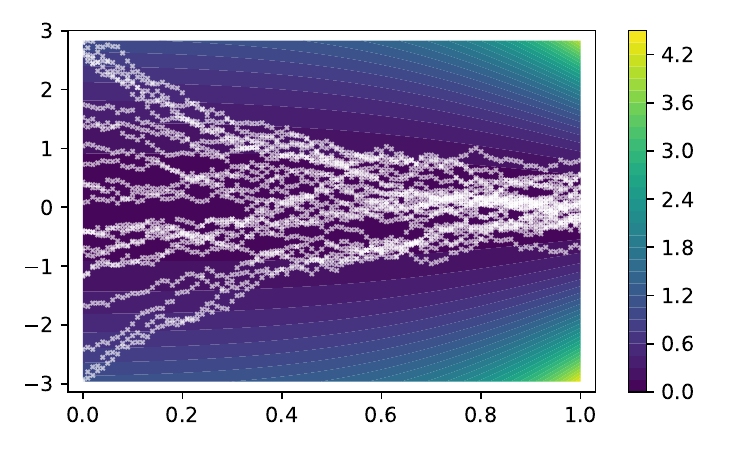}
    \includegraphics[width=0.4\textwidth]{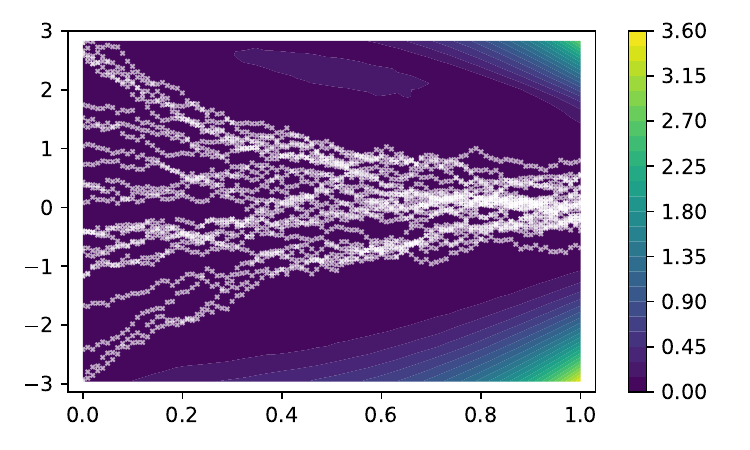}
    \caption{
    Numerical results of \Cref{alg_mini} for the interbank systemic-risk MFG in \cref{sec_sr}. Unless otherwise noted, curves and shaded bands are the mean and mean $\pm 2\times\mathrm{SD}$ over five independent runs.
    (Top left) $\mathrm{RE}_{\infty}$ and $\mathrm{RE}_1$ versus iteration.
    (Top middle) $\mathrm{RC}$ versus iteration.
    (Top right) Empirical mean reserve $\hat m_t$
    and the exact counterpart $m_t^*$, with horizontal axis $t$ and vertical axis mean reserve.
    (Middle) Pointwise mean learned controls at $t=0$, $0.5$, and $0.99$; the rug marks indicate the aggregated empirical distribution of the $5 \times 2049$ simulated reserve states $Z^{\mu_I, u_{\alpha}}$; the horizontal and vertical axes represent $z$ and control value, respectively.
    (Bottom) For one representative run, the left and right panels show the predicted value landscape $v_{\theta}$ and the absolute error $\abs{v_{\theta}-v}$, with simulated paths $Z^{\mu_I, u_{\alpha}}$ overlaid; the horizontal and vertical axes represent $t$ and $z$, respectively.
    }
    \label{fig_sr_results}
\end{figure}

As a benchmark for numerical performance, we present the explicit solution to this model.
The value function $v$ and the optimal control $u^*$ for this model are given by
\begin{equation}\label{eq_sr_ansatz}
    v(x) = \frac{1}{2}\eta(t)\bigl(z - \bar{m}_t\bigr)^2 + \gamma(t), \quad u^*(x) = \bigl(c_3 + \eta(t)\bigr)\bigl(\bar{m}_t - z\bigr)
\end{equation}
for $x = (t, z) \in [0, T) \times \R$, where the functions $\eta$ and $\gamma$ satisfy the following system of ODEs:
\begin{equation}\label{eq_sr_riccati}
    \dot{\eta}(t) = 2(c_1 + c_3)\eta(t) + \eta(t)^2 - (c_4 - c_3^2), \quad \dot{\gamma}(t) = -\frac{c_2^2}{2}\eta(t), \quad \eta(T) = c_5, \quad \gamma(T) = 0
\end{equation}
with $t \in [0, T)$.
As in \Cref{sec_numerical_lq}, the above ODEs are solved numerically using the same variable-step Runge--Kutta solver settings to obtain reference solutions for $v$ and $u^*$.
The results of \Cref{alg_mini} for this example are shown in \Cref{fig_sr_results}.

The numerical results in \Cref{fig_sr_results} show that \Cref{alg_mini} performs well for this systemic-risk example.
As seen in the top-left and top-middle panels, the mean values of $\mathrm{RE}_1$, $\mathrm{RE}_{\infty}$, and $\mathrm{RC}$ decrease steadily over the iterations and, by the final iteration, reach $1.48 \times 10^{-2}$, $2.64 \times 10^{-2}$, and $1.32 \times 10^{-2}$, respectively, with small standard deviations across five independent runs.
The top-right panel presents a comparison between the empirical mean reserve $\hat m_t$ under the learned equilibrium measure and the exact mean reserve trajectory $m_t^*$.
The close agreement between $\hat m_t$ and $m_t^*$ indicates that the learned equilibrium measure accurately captures the population distribution of the state process under the optimal control.

The middle panel compares the learned and exact controls, with rug marks showing the aggregated empirical distribution of the simulated reserve states $Z^{\mu_I, u_{\alpha}}$.
The learned control $u_{\alpha}(t, z)$ closely matches the exact optimal control $u^*(t, z)$ for $z$ in the region where the simulated paths concentrate, while larger discrepancies appear in the tails of the distribution of $Z_t^{\mu_I, u_{\alpha}}$, which are sampled less often.
The state distribution also becomes more concentrated as $t$ approaches the terminal time $T$. 
This is consistent with the cost functional in \eqref{eq_sr_cost}, which encourages agents to stay close to the population mean $\bar{m}_t$ and therefore leads to a tighter distribution of $Z_t^{\mu_I, u_{\alpha}}$ for $t$ near $T$, as shown in the bottom panels of \Cref{fig_sr_results}.

\subsection{Target-tracking MFG with non-linear interactions}\label{sec_tar_mfg}

\begin{figure}[t]
    \centering
    \includegraphics[width=0.42\textwidth]{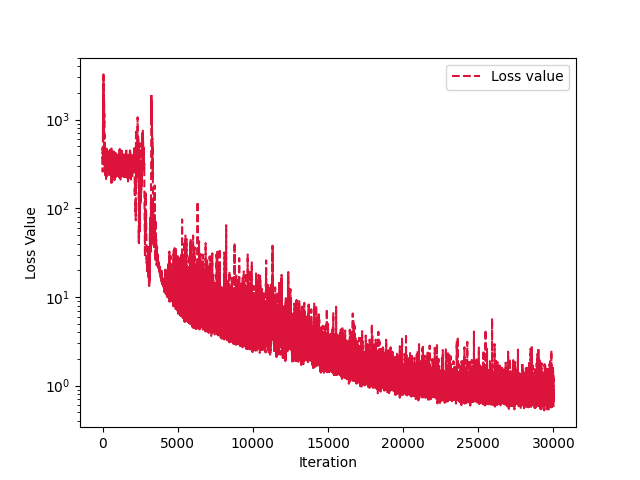}
    \includegraphics[width=0.38\textwidth]{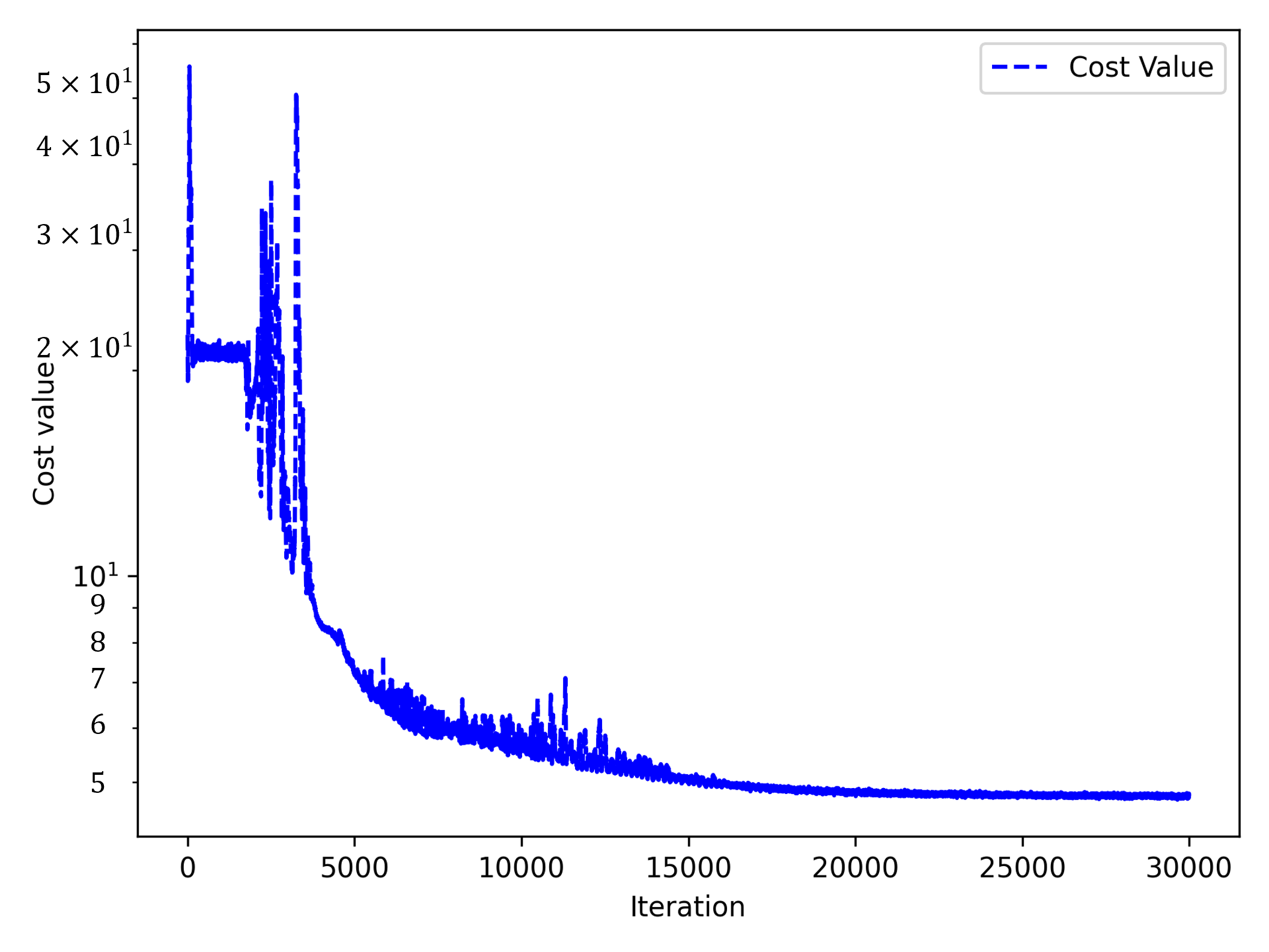}\\
	\includegraphics[width=0.43\textwidth]{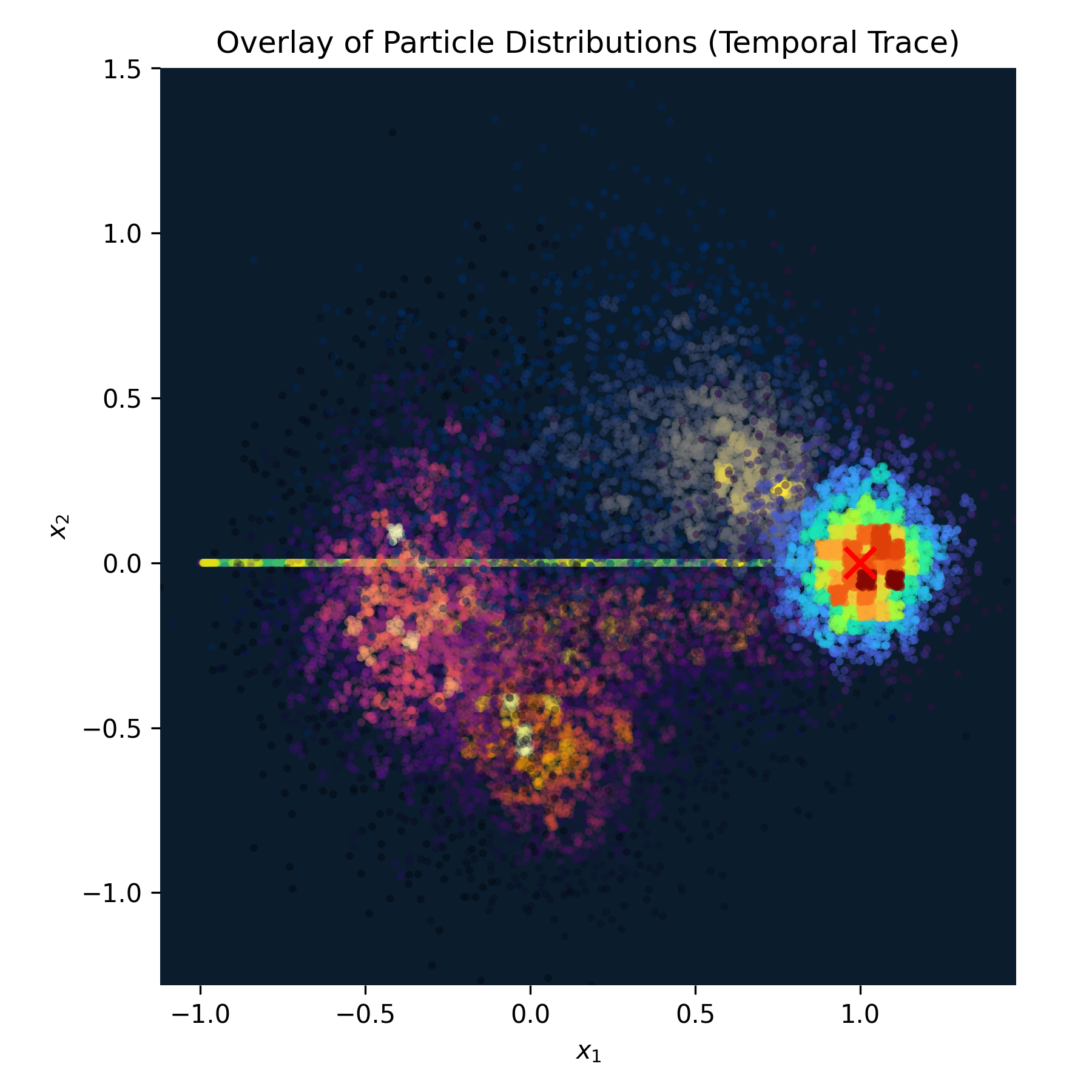}
	\includegraphics[width=0.54\textwidth]{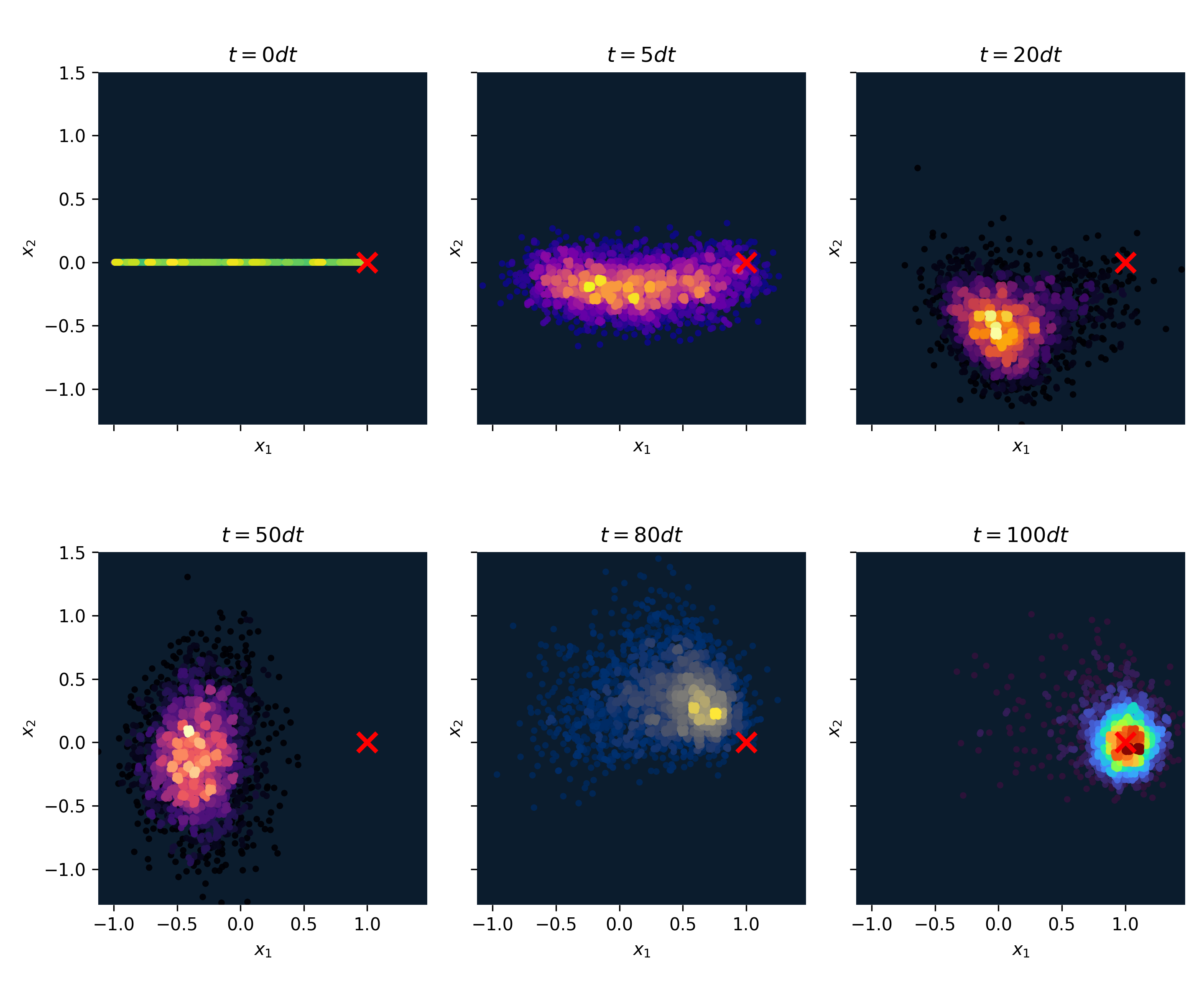}
    \caption{
    Results of \Cref{alg_mini} for the MFG in \Cref{sec_tar_mfg}.
    (Upper left) Loss versus iteration.
    (Upper right) Empirical value $\hat{J}$ of the cost functional in \eqref{eq_cost} versus iteration.
    (Lower left) Overlay of snapshots of agent positions at $t = i/5 T$, $i = 0, 1, \cdots, 5$.
    (Lower right) Agent trajectories over time, with snapshots at $t = i/5 T$, $i = 0, 1, \cdots, 5$; the red ``x'' marks the target position $z^*(T)$.}
    \label{fig_tar}
\end{figure}

This example considers a target-tracking MFG with non-linear interactions and thus serves as a test case for the proposed method in a more general setting.
Specifically, the state dynamics and cost functional are given by \eqref{eq_state} and \eqref{eq_cost}, respectively, with the following specifications:
\begin{align}
    &b_z(x, \mu, u) = u, \quad \sigma_z(x, \mu, u) = 0.5 I_d, \quad g(x, \mu) = 0.5 c_4 |z - z^*(T)|^2, \label{eq_tar_bsg}\\
    &f(x, \mu, u) = c_1 \abs{z - z^*(t)}^2 + c_2 \abs{u}^2 + F(x, \mu),  \quad T = 1, \label{eq_tar_f}
\end{align}
for $x = (t, z) \in [0, T] \times \R^d$ and $u \in \R^d$. 
The definitions of $f$ and $g$ encourage agents to track the target trajectory $z^*: [0, T] \to \R^d$ which is a circular path defined by
\begin{equation*}
    z_i^*(t) := \sin(2 \pi t + i \pi/2), \quad i = 1, 2, \cdots, d, \quad t \in [0, T]; 
\end{equation*}
The interaction term $F(x, \mu)$ is introduced to prevent agents from aggregating at the same location, and is defined as
\begin{equation}\label{eq_defFxm}
    F(x, \mu) := \int_{[0, T] \times \R^d} \delta_t(s) \exp\big(\!- c_F \abs{z - y}^2\big) \mu(\!\di s \times \di y), \quad x = (t, z) \in [0, T] \times \R^d,
\end{equation}
We set $c_1 = 1$, $c_2 = 0.1$, $c_4 = 10$, $c_F = 1$.
For visualization, we set the spatial dimension to $d = 2$.
The initial state is given by $X_0 = (0, Z_0)$, where $Z_0$ is uniformly distributed on the line segment $\{(s, 0) : s \in [-1, 1]\}$.
The results of \Cref{alg_mini} for this example are shown in \Cref{fig_tar}.

As shown in the upper-right subplot of \Cref{fig_tar}, the empirical cost decreases by about $0.1$ from the initial to the final iteration, indicating that the learned control $u_{\alpha}$ substantially improves upon the initial guess.
The lower two subplots of \Cref{fig_tar} show that the agents move in an approximately circular pattern over time, consistent with the circular target trajectory, and remain close to the target at the terminal time $T$.
These observations suggest that the proposed method performs reasonably well in this more complex setting, although the absence of an explicit solution precludes a more quantitative assessment of the accuracy of the learned control and value function.

\subsection{MFG with nonlinear state dynamics and interactions}\label{sec_numerical_barrier}

\begin{figure}[t]
    \centering
	\includegraphics[width=0.4\textwidth]{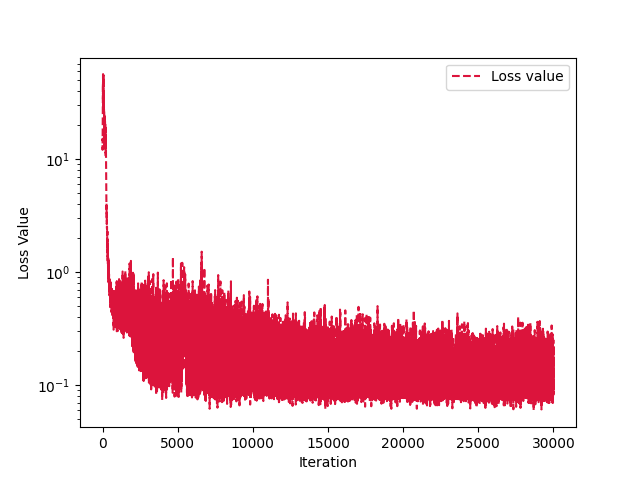}
	\includegraphics[width=0.4\textwidth]{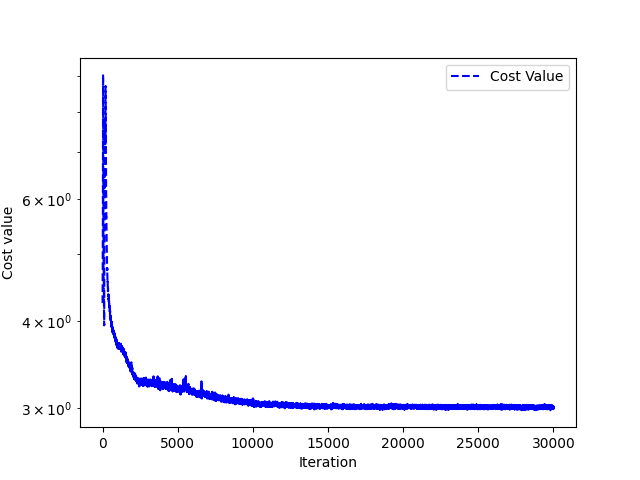}
	\includegraphics[width=0.4\textwidth]{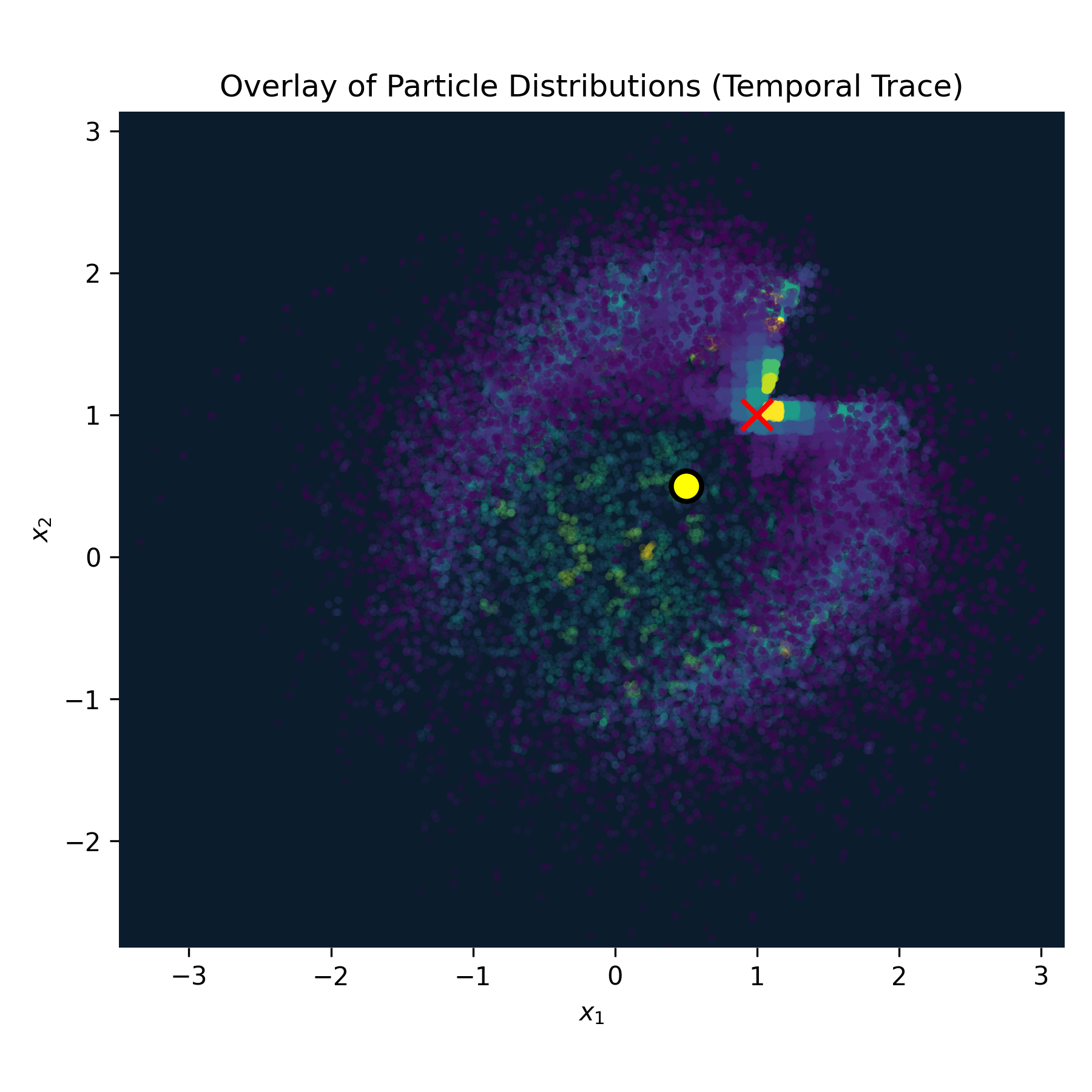}
	\includegraphics[width=0.5\textwidth]{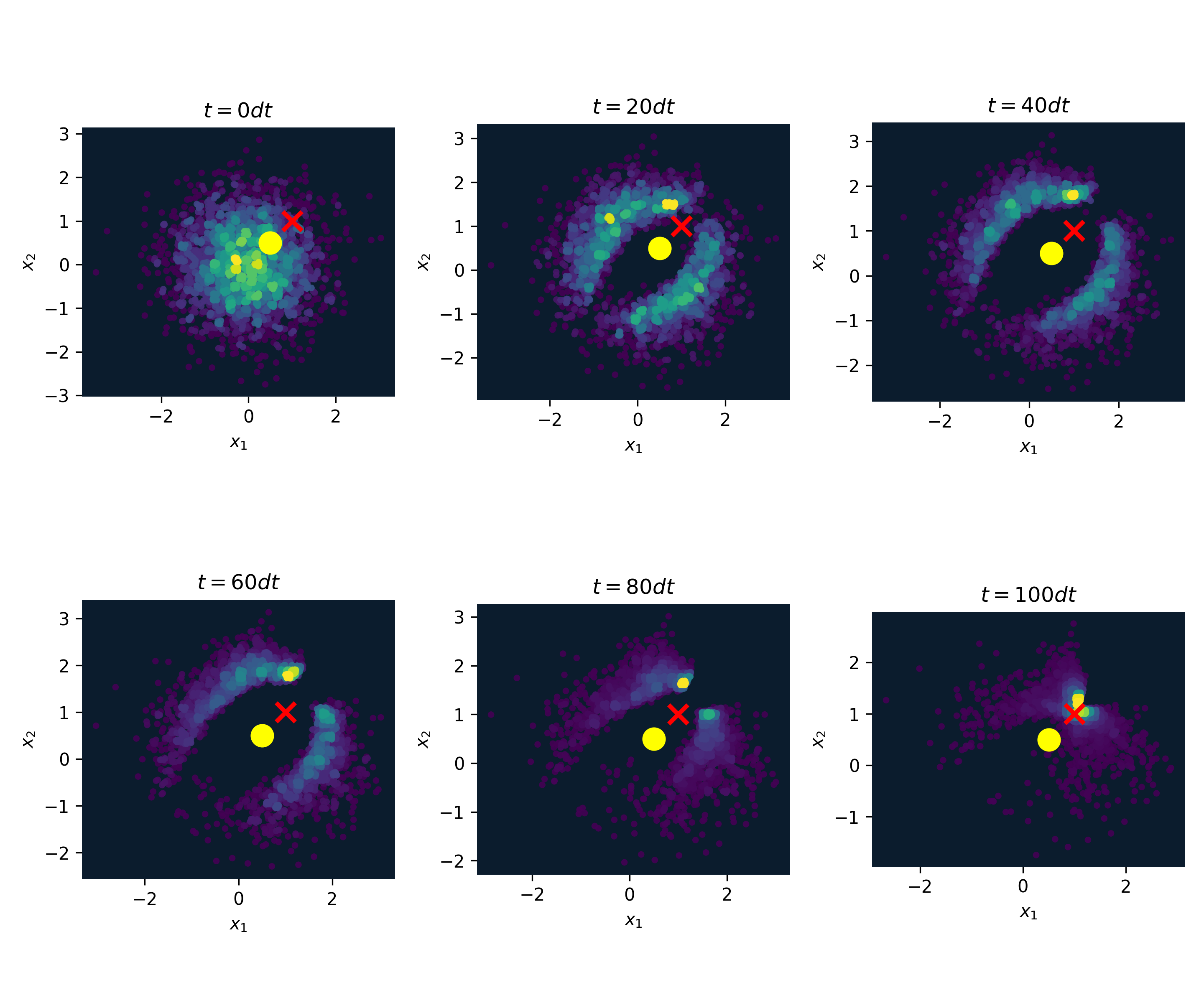}
    \caption{
    Results of \Cref{alg_mini} for the MFG in \Cref{sec_numerical_barrier}.
    (Upper left) Loss versus iteration.
    (Upper right) Empirical value $\hat{J}$ of the cost functional in \eqref{eq_cost} versus iteration.
    (Lower left) Overlay of snapshots of agent positions at $t = i/5 T$, $i = 0, 1, \cdots, 5$.
    (Lower right) Agent trajectories with snapshots at $t = i/5 T$, $i = 0, 1, \cdots, 5$; the red ``x'' marks the target position $z^*(T)$, and the yellow ``o'' denotes the barrier point.
    }
    \label{fig_bar}
\end{figure}

This example considers an MFG with nonlinear state dynamics and interactions to demonstrate the applicability of the proposed method in a more complex setting.
The state dynamics and the cost functional are given by \eqref{eq_state} and \eqref{eq_cost}, respectively, with the following specifications:
\begin{align}
    &b_z(x, \mu, u) = u, \quad \sigma_z(x, \mu, u) = \frac{c_1}{\sqrt{d}} \sum_{i=1}^d \sin(i + z_i) I_d, \quad g(x, \mu) = 0.5 c_4 |z - z^*(T)|^2, \label{eq_barrier_bsg}\\
    &f(x, \mu, u) = c_0 / (1 + c_2 |z - 0.5\B{1}_d|^2) + c_3 |u|^2 + F(x, \mu),  \quad T = 1 \label{eq_barrier_f}
\end{align}
for $x = (t, z) \in [0, T] \times \R^d$ and $u \in \R^d$, where the interaction term $F(x, \mu)$ is defined in \eqref{eq_defFxm} with $c_F = 1$.
The constants are set to $c_1 = c_3 = 0.1$, $c_2 = c_4 = 1$, and $c_0 = 5$.
In \eqref{eq_barrier_f}, the first term in $f$ penalizes states near the barrier point $0.5 \B{1}_d$, thereby encouraging agents to avoid it, while the third term, as in \Cref{sec_tar_mfg}, discourages aggregation at the same location.
For visualization, we set the spatial dimension to $d = 2$.
The initial state is given by $X_0 = (0, Z_0)$, where $Z_0$ is a two-dimensional standard Gaussian random variable.
The results of \Cref{alg_mini} for this example are shown in \Cref{fig_bar}.

As shown in the lower subplots of \Cref{fig_bar}, the agents initially spread according to the Gaussian distribution of the initial state, then move outward to avoid the barrier point $0.5 \B{1}_d$, where the running cost is high, and finally approach the target position $z^*(T)$ at the terminal time.
This behavior is consistent with the specified state dynamics and cost functional.
Together with the example in \Cref{sec_tar_mfg}, this example illustrates the effectiveness of the proposed method in handling MFG with different forms of nonlinearity and interaction.

\section{Conclusion}\label{sec_conclusion}

In this paper, we proposed a deep policy iteration method for high-dimensional finite-horizon mean-field games based on a regenerative reformulation. By recasting the original problem as a regenerative MFG with deterministic cycles and a resetting mechanism, we obtained a cycle-wise scheme consisting of policy evaluation, policy improvement, and population measure estimation. The equilibrium measure is approximated by an empirical particle measure updated through a one-step random mapping induced by the Euler--Maruyama discretization, thereby avoiding full trajectory simulation over the entire time horizon at each iteration. The discretized PE and PI steps are written in weak form with neural networks representing the value function and feedback control and an adversarial test-function network providing a Galerkin-type formulation. This avoids directly solving the coupled HJB--FP system, explicitly computing conditional expectations in policy evaluation, evaluating higher-order derivatives of the HJB equation, and performing pointwise Hamiltonian minimization in policy improvement.

We tested the method on linear-quadratic, interbank systemic-risk, target-tracking, and nonlinear MFG examples. The results indicate strong performance in both low- and high-dimensional settings and computational feasibility up to dimension 10{,}000. Future work includes establishing rigorous convergence and error estimates, extending the framework to more general mean-field interactions and common noise, and developing adaptive sampling and network architectures tailored to specific problem structures.

\appendix
\section{Derivation of the explicit solutions in \eqref{eq_ansatz_lq}--\eqref{eq_final_ode_g}}\label{appendix_lq}

We now provide a detailed derivation of the explicit solution \eqref{eq_ansatz_lq}--\eqref{eq_final_ode_g}.
For the linear-quadratic instance of \eqref{eq_state} and \eqref{eq_cost} with specifications \eqref{eq_state_lq} and \eqref{eq_fg_lg}, the value function $v$ satisfies the following HJB equation \cite[p. 268, Theorem 5.1]{Yong1999Stochastic}:
\begin{equation}\label{eq_hjb_lq}
\begin{aligned}
    \partial_t v + \inf_{\kappa \in \R^d} \big\{ \kappa^{\top} \partial_z v + \frac{c_2}{2} \abs{\kappa}^2 \big\} + \big\{c_0 (\bar{m}_t - z) + c_1 \bigl(z^*(t) - \bar{m}_t\bigr) \big\}^{\top} \partial_z v &\\
    + \frac{c_{\sigma}^2}{2} \Delta_z v + \frac{c_3}{2} \abs{z - \bar{m}_t}^2 + \frac{c_4}{2} \abs{z - z^*(t)}^2 &= 0
\end{aligned}
\end{equation}
for $(t, z) \in [0, T) \times \R^d$, with terminal condition $v(T, z) = c_5/2 \abs{z - z^*(T)}^2 + 1/2$ for $z \in \R^d$, where $\Delta_z v$ denotes the Laplacian of $v$ with respect to $z$.
The optimal feedback control is given by
\begin{equation}\label{eq_opt_lq}
    u^*(t, z) = - c_2^{-1} \partial_z v(t, z), \quad (t, z) \in [0, T] \times \R^d.
\end{equation}
Substituting \eqref{eq_opt_lq} into the governing equation in \eqref{eq_hjb_lq}, we obtain
\begin{equation}\label{eq2_hjb_lp}
\begin{aligned}
    \partial_t v + \big\{c_0 (\bar{m}_t - z) + c_1 \bigl(z^*(t) - \bar{m}_t\bigr) \big\}^{\top} \partial_z v - \frac{1}{2 c_2} \abs{\partial_z v}^2& \\
    + \frac{c_{\sigma}^2}{2} \Delta_z v + \frac{c_3}{2} \abs{z - \bar{m}_t}^2 + \frac{c_4}{2} \abs{z - z^*(t)}^2 &= 0
\end{aligned}
\end{equation}
for $(t, z) \in [0, T) \times \R^d$.

Having established that the value function $v$ admits the form given in \eqref{eq_ansatz_lq}, the optimal control $u^*$ in \eqref{eq_opt_lq} follows directly from this ansatz. 
We therefore seek a quadratic solution of the form \eqref{eq_ansatz_lq}, where $a: [0, T] \to \R$, $b: [0, T] \to \R^d$, and $\gamma: [0, T] \to \R$ are functions to be determined.
To this end, with the ansatz \eqref{eq_ansatz_lq}, we compute $\partial_t v$, $\partial_z v$, and $\Delta_z v$ as follows:
\begin{align}
    &\partial_t v(t, z) = \frac{1}{2} a'(t) \abs{z - \bar{m}_t}^2 - a(t) (z - \bar{m}_t)^{\top} \dot{\bar{m}}_t + b'(t)^{\top} (z - \bar{m}_t)- b(t)^{\top} \dot{\bar{m}}_t + \gamma'(t), \label{eq_ptv}\\
    &\partial_z v(t, z) = a(t)(z - \bar{m}_t) + b(t), \qquad \Delta_z v(t, z) = d\, a(t). \label{eq_dxv}
\end{align}
Inserting the above expressions into \eqref{eq2_hjb_lp} gives
\begin{equation}\label{eq_expanded_hjb_lq}
\begin{aligned}
    \frac{1}{2} a'(t) \abs{z - \bar{m}_t}^2 - a(t) (z - \bar{m}_t)^{\top} \dot{\bar{m}}_t + b'(t)^{\top} (z - \bar{m}_t)- b(t)^{\top} \dot{\bar{m}}_t + \gamma'(t) &\\
    + \big\{c_0 (\bar{m}_t - z) + c_1 \bigl(z^*(t) - \bar{m}_t\bigr) \big\}^{\top} \big\{a(t)(z - \bar{m}_t) + b(t)\big\} &  \\
     - \frac{1}{2 c_2} \abs{a(t)(z - \bar{m}_t) + b(t)}^2 + \frac{c_{\sigma}^2}{2} d\, a(t) + \frac{c_3}{2} \abs{z - \bar{m}_t}^2 + \frac{c_4}{2} \abs{z - z^*(t)}^2 &= 0,
\end{aligned}
\end{equation}
where the last term on the left-hand side can be rewritten as
\begin{equation*}
    \frac{c_4}{2} \abs{z - z^*(t)}^2 = \frac{c_4}{2} \abs{z - \bar{m}_t}^2 + c_4 (z - \bar{m}_t)^{\top} \bigl(\bar{m}_t - z^*(t)\bigr) + \frac{c_4}{2} \abs{\bar{m}_t - z^*(t)}^2.
\end{equation*}

Substituting the above expression into \eqref{eq_expanded_hjb_lq} and rearranging terms by collecting the coefficients of $\abs{z - \bar{m}_t}^2$, $z - \bar{m}_t$, and the constant term yields
\begin{align}
    \frac{1}{2} a'(t) - c_0 a(t) - \frac{1}{2 c_2} a(t)^2 + \frac{c_3}{2} + \frac{c_4}{2} &= 0, \label{eq_coef_quad} \\
    -a(t)\dot{\bar{m}}_t + b'(t) - c_0 b(t) + a(t)c_1\bigl(z^*(t)-\bar{m}_t\bigr) - c_2^{-1} a(t)b(t) + c_4\bigl(\bar{m}_t - z^*(t)\bigr) &= 0, \label{eq_coef_lin} \\
    -b(t)^\top\dot{\bar{m}}_t + \gamma'(t) + c_1\bigl(z^*(t)-\bar{m}_t\bigr)^\top b(t) - \frac{1}{2 c_2}|b(t)|^2 + \frac{c_{\sigma}^2}{2}d\,a(t) + \frac{c_4}{2}|\bar{m}_t - z^*(t)|^2 &= 0. \label{eq_coef_const}
\end{align}
To obtain the terminal conditions for $a$, $b$, and $\gamma$, 
we rewrite the terminal condition $v(T, z) = c_5/2 \abs{z - z^*(T)}^2 + 1/2$ as 
\begin{equation*}
    v(T, z) = \frac{c_5}{2} \abs{z - \bar{m}_T}^2 + c_5 \bigl(\bar{m}_T - z^*(T)\bigr)^{\top} (z - \bar{m}_T) + \frac{c_5}{2} \abs{\bar{m}_T - z^*(T)}^2 + \frac{1}{2}, \quad z \in \R^d.
\end{equation*}
Recalling the ansatz \eqref{eq_ansatz_lq} and matching the coefficients of $\abs{z - \bar{m}_T}^2$, $z - \bar{m}_T$, and the constant term in the above equation, 
we obtain the terminal conditions for $a$, $b$, and $\gamma$ as follows:
\begin{equation}\label{eq_term}
    a(T) = c_5, \quad b(T) = c_5 \bigl(\bar{m}_T - z^*(T)\bigr), \quad \gamma(T) = \frac{c_5}{2} \abs{z^*(T) - \bar{m}_T}^2.
\end{equation}

The ODE~\eqref{eq_final_ode_a} for $a$ is obtained by simplifying \eqref{eq_coef_quad}, with the terminal condition specified in \eqref{eq_term}.

To simplify \eqref{eq_coef_lin}, we use the expectation of the optimal state dynamics.
Under the optimal feedback control $u^*$ given in the second equation in \eqref{eq_ansatz_lq}, the state equation becomes
\begin{equation*}
    Z_t = Z_0 + \int_0^t \bigl\{(-c_0 - a(s))(Z_s - \bar{m}_s) - b(s) + c_1\bigl(z^*(s)-\bar{m}_s\bigr)\bigr\} \di s + \int_0^t c_{\sigma} I_d \, \di B_s, \quad t \in [0, T].
\end{equation*}
Taking expectation on both sides of the above equation, using $\mathbb{E}[Z_t - \bar{m}_t] = 0$ to simplify the integrand, and then differentiating with respect to $t$, we obtain \eqref{eq_final_ode_m} for $\bar{m}$ with the initial condition specified by the definition $\bar{m}_t := \E{Z_t}$, $t \in [0, T]$.

Using $\dot{\bar{m}}_t + c_2^{-1} b(t) = c_1\bigl(z^*(t)-\bar{m}_t\bigr)$ implied by \eqref{eq_final_ode_m} to simplify \eqref{eq_coef_lin} yields the ODE \eqref{eq_final_ode_b} for $b(t)$, where the terminal condition is from the second condition in \eqref{eq_term}.

Finally, substituting $\dot{\bar{m}}_t$ from \eqref{eq_final_ode_m} into \eqref{eq_coef_const}, and using the relation $c_2^{-1} b(t) = -\dot{\bar{m}}_t + c_1\bigl(z^*(t)-\bar{m}_t\bigr)$ (also implied by \eqref{eq_final_ode_m}), we obtain the ODE~\eqref{eq_final_ode_g} for $\gamma$, where the terminal condition is from the third condition in \eqref{eq_term}.

\bibliographystyle{plain}
\bibliography{bibliography}

\begin{thebibliography}{10}

\bibitem{Achdou19}
Yves Achdou and Jean-Michel Lasry.
\newblock Mean field games for modeling crowd motion.
\newblock In {\em Contributions to partial differential equations and applications}, volume~47 of {\em Comput. Methods Appl. Sci.}, pages 17--42. Springer, Cham, 2019.

\bibitem{Al2022Extensions}
Ali Al-Aradi, Adolfo Correia, Gabriel Jardim, Danilo de~Freitas Naiff, and Yuri Saporito.
\newblock Extensions of the deep {G}alerkin method.
\newblock {\em Appl. Math. Comput.}, 430:Paper No. 127287, 18, 2022.

\bibitem{Andersson}
Daniel Andersson and Boualem Djehiche.
\newblock A maximum principle for {SDE}s of mean-field type.
\newblock {\em Appl. Math. Optim.}, 63(3):341--356, 2011.

\bibitem{Assouli2024Deep}
Mouhcine Assouli and Badr Missaoui.
\newblock Deep policy iteration for high-dimensional mean field games.
\newblock {\em Appl. Math. Comput.}, 481:Paper No. 128923, 12, 2024.

\bibitem{Bensoussan13}
Alain Bensoussan, Jens Frehse, and Phillip Yam.
\newblock {\em Mean field games and mean field type control theory}.
\newblock SpringerBriefs in Mathematics. Springer, New York, 2013.

\bibitem{Bensoussan18Mean}
Alain Bensoussan, Tao Huang, and Mathieu Lauri\`ere.
\newblock Mean field control and mean field game models with several populations.
\newblock {\em Minimax Theory Appl.}, 3(2):173--209, 2018.

\bibitem{cai2024martingale}
Wei Cai, Shuixin Fang, Wenzhong Zhang, and Tao Zhou.
\newblock Martingale deep learning for very high dimensional quasi-linear partial differential equations and stochastic optimal controls.
\newblock {\em arXiv preprint arXiv:2408.14395}, 2024.

\bibitem{cai2024socmartnet}
Wei Cai, Shuixin Fang, and Tao Zhou.
\newblock {SOC-MartNet}: A martingale neural network for the hamilton-jacobi-bellman equation without explicit $\inf_{u \in U} {H}$ in stochastic optimal controls.
\newblock {\em SIAM J. Sci. Comput.}, 47(4):C795--C819, 2025.

\bibitem{Cai2026Deep}
Wei Cai, Shuixin Fang, and Tao Zhou.
\newblock Deep random difference method for high-dimensional quasilinear parabolic partial differential equations.
\newblock {\em J. Comput. Phys.}, page 114767, 2026.

\bibitem{Andrew26}
Wei Cai, Andrew He, and Daniel Margolis.
\newblock Deep{M}art{N}et: a {M}artingale-based deep neural network learning method for {D}irichlet {BVP}s and eigenvalue problems of elliptic {PDE}s in $\mathbb{R}^d$.
\newblock {\em SIAM J. Sci. Comput.}, 48(1):C25--C50, 2026.

\bibitem{Lasry13}
P.~Cardaliaguet, J.-M. Lasry, P.-L. Lions, and A.~Porretta.
\newblock Long time average of mean field games with a nonlocal coupling.
\newblock {\em SIAM J. Control Optim.}, 51(5):3558--3591, 2013.

\bibitem{cardaliaguet2010notes}
Pierre Cardaliaguet.
\newblock Notes on mean field games.
\newblock Technical report, Technical report Technical report, 2010.

\bibitem{Cardaliaguet2017Learning}
Pierre Cardaliaguet and Saeed Hadikhanloo.
\newblock Learning in mean field games: the fictitious play.
\newblock {\em ESAIM Control Optim. Calc. Var.}, 23(2):569--591, 2017.

\bibitem{Carmona18}
Ren\'{e} Carmona and Fran\c{c}ois Delarue.
\newblock {\em Probabilistic theory of mean field games with applications. {I}}, volume~83 of {\em Probability Theory and Stochastic Modelling}.
\newblock Springer, Cham, 2018.
\newblock Mean field FBSDEs, control, and games.

\bibitem{Carmona15Mean}
Ren\'{e} Carmona, Jean-Pierre Fouque, and Li-Hsien Sun.
\newblock Mean field games and systemic risk.
\newblock {\em Commun. Math. Sci.}, 13(4):911--933, 2015.

\bibitem{Carmona15probabilistic}
Ren\'{e} Carmona and Daniel Lacker.
\newblock A probabilistic weak formulation of mean field games and applications.
\newblock {\em Ann. Appl. Probab.}, 25(3):1189--1231, 2015.

\bibitem{Elliott13}
Robert Elliott, Xun Li, and Yuan-Hua Ni.
\newblock Discrete time mean-field stochastic linear-quadratic optimal control problems.
\newblock {\em Automatica J. IFAC}, 49(11):3222--3233, 2013.

\bibitem{Gao2023Failure}
Zhiwei Gao, Liang Yan, and Tao Zhou.
\newblock Failure-informed adaptive sampling for {PINN}s.
\newblock {\em SIAM J. Sci. Comput.}, 45(4):A1971--A1994, 2023.

\bibitem{Garnier13}
Josselin Garnier, George Papanicolaou, and Tzu-Wei Yang.
\newblock Large deviations for a mean field model of systemic risk.
\newblock {\em SIAM J. Financial Math.}, 4(1):151--184, 2013.

\bibitem{germain2022Approximation}
Maximilien Germain, Huy\^{e}n Pham, and Xavier Warin.
\newblock Approximation error analysis of some deep backward schemes for nonlinear {PDE}s.
\newblock {\em SIAM J. Sci. Comput.}, 44(1):A28--A56, 2022.

\bibitem{Han2020Deep}
Jiequn Han and Ruimeng Hu.
\newblock Deep fictitious play for finding {M}arkovian {N}ash equilibrium in multi-agent games.
\newblock In Jianfeng Lu and Rachel Ward, editors, {\em Proceedings of The First Mathematical and Scientific Machine Learning Conference}, volume 107 of {\em Proceedings of Machine Learning Research}, pages 221--245. PMLR, 20--24 Jul 2020.

\bibitem{Han2024Learning}
Jiequn Han, Ruimeng Hu, and Jihao Long.
\newblock Learning high-dimensional {M}c{K}ean-{V}lasov forward-backward stochastic differential equations with general distribution dependence.
\newblock {\em SIAM J. Numer. Anal.}, 62(1):1--24, 2024.

\bibitem{han2018solving}
Jiequn Han, Arnulf Jentzen, and Weinan E.
\newblock Solving high-dimensional partial differential equations using deep learning.
\newblock {\em Proceedings of the National Academy of Sciences}, 115(34):8505--8510, 2018.

\bibitem{He2023Learning}
Di~He, Shanda Li, Wenlei Shi, Xiaotian Gao, Jia Zhang, Jiang Bian, Liwei Wang, and Tie-Yan Liu.
\newblock Learning physics-informed neural networks without stacked back-propagation.
\newblock In Francisco Ruiz, Jennifer Dy, and Jan-Willem van~de Meent, editors, {\em Proceedings of The 26th International Conference on Artificial Intelligence and Statistics}, volume 206 of {\em Proceedings of Machine Learning Research}, pages 3034--3047. PMLR, 25--27 Apr 2023.

\bibitem{Hu2021Deep}
Ruimeng Hu.
\newblock Deep fictitious play for stochastic differential games.
\newblock {\em Commun. Math. Sci.}, 19(2):325--353, 2021.

\bibitem{Hu2024Hutchinson}
Zheyuan Hu, Zekun Shi, George~Em Karniadakis, and Kenji Kawaguchi.
\newblock Hutchinson trace estimation for high-dimensional and high-order physics-informed neural networks.
\newblock {\em Comput. Methods Appl. Mech. Engrg.}, 424:Paper No. 116883, 17, 2024.

\bibitem{hu2024sdgd}
Zheyuan Hu, Khemraj Shukla, George~Em Karniadakis, and Kenji Kawaguchi.
\newblock Tackling the curse of dimensionality with physics-informed neural networks.
\newblock {\em Neural Networks}, 176:106369, 2024.

\bibitem{Hu2025Bias}
Zheyuan Hu, Zhouhao Yang, Yezhen Wang, George~E. Karniadakis, and Kenji Kawaguchi.
\newblock Bias-{V}ariance {T}rade-{O}ff in {P}hysics-{I}nformed {N}eural {N}etworks with {R}andomized {S}moothing for {H}igh-{D}imensional {PDE}s.
\newblock {\em SIAM J. Sci. Comput.}, 47(4):C846--C872, 2025.

\bibitem{Minyi10}
Minyi Huang.
\newblock Large-population {LQG} games involving a major player: the {N}ash certainty equivalence principle.
\newblock {\em SIAM J. Control Optim.}, 48(5):3318--3353, 2009/10.

\bibitem{Minyi12}
Minyi Huang, Peter~E. Caines, and Roland~P. Malham\'{e}.
\newblock Social optima in mean field {LQG} control: centralized and decentralized strategies.
\newblock {\em IEEE Trans. Automat. Control}, 57(7):1736--1751, 2012.

\bibitem{hure2020deep}
C\^{o}me Hur\'{e}, Huy\^{e}n Pham, and Xavier Warin.
\newblock Deep backward schemes for high-dimensional nonlinear {PDE}s.
\newblock {\em Math. Comp.}, 89(324):1547--1579, 2020.

\bibitem{Yanwei2022Policy}
Yanwei Jia and Xun~Yu Zhou.
\newblock Policy evaluation and temporal-difference learning in continuous time and space: A martingale approach.
\newblock {\em Journal of Machine Learning Research}, 23(154):1--55, 2022.

\bibitem{Yanwei2023Policy}
Yanwei Jia and Xun~Yu Zhou.
\newblock Policy gradient and actor-critic learning in continuous time and space: Theory and algorithms.
\newblock {\em Journal of Machine Learning Research}, 23(275):1--50, 2022.

\bibitem{Achim2020Probability}
Achim Klenke.
\newblock {\em Probability Theory}.
\newblock Springer Cham, third edition, 2020.

\bibitem{Kloeden1992Numerical}
Peter~E. Kloeden and Eckhard Platen.
\newblock {\em Numerical solution of stochastic differential equations}, volume~23 of {\em Applications of Mathematics (New York)}.
\newblock Springer-Verlag, Berlin, 1992.

\bibitem{Lachapelle16}
Aim\'{e} Lachapelle, Jean-Michel Lasry, Charles-Albert Lehalle, and Pierre-Louis Lions.
\newblock Efficiency of the price formation process in presence of high frequency participants: a mean field game analysis.
\newblock {\em Math. Financ. Econ.}, 10(3):223--262, 2016.

\bibitem{Lachapelle10}
Aime Lachapelle, Julien Salomon, and Gabriel Turinici.
\newblock Computation of mean field equilibria in economics.
\newblock {\em Math. Models Methods Appl. Sci.}, 20(4):567--588, 2010.

\bibitem{lachapelle2011mean}
Aim{\'e} Lachapelle and Marie-Therese Wolfram.
\newblock On a mean field game approach modeling congestion and aversion in pedestrian crowds.
\newblock {\em Transportation research part B: methodological}, 45(10):1572--1589, 2011.

\bibitem{Li2024Neural}
Xingjian Li, Deepanshu Verma, and Lars Ruthotto.
\newblock A neural network approach for stochastic optimal control.
\newblock {\em SIAM J. Sci. Comput.}, 46(5):C535--C556, 2024.

\bibitem{Liang2024Actor}
Hong Liang, Zhiping Chen, and Kaili Jing.
\newblock Actor-critic reinforcement learning algorithms for mean field games in continuous time, state and action spaces.
\newblock {\em Appl. Math. Optim.}, 89(3):Paper No. 72, 35, 2024.

\bibitem{Lin2021Alternating}
Alex~Tong Lin, Samy~Wu Fung, Wuchen Li, Levon Nurbekyan, and Stanley~J. Osher.
\newblock Alternating the population and control neural networks to solve high-dimensional stochastic mean-field games.
\newblock {\em Proc. Natl. Acad. Sci. USA}, 118(31):Paper No. e2024713118, 10, 2021.

\bibitem{Liu2020Multi}
Ziqi Liu, Wei Cai, and Zhi-Qin~John Xu.
\newblock Multi-scale deep neural network ({M}scale{DNN}) for solving {P}oisson-{B}oltzmann equation in complex domains.
\newblock {\em Commun. Comput. Phys.}, 28(5):1970--2001, 2020.

\bibitem{Mou2024Bellman}
Wenlong Mou and Yuhua Zhu.
\newblock On bellman equations for continuous-time policy evaluation i: discretization and approximation, 2024.

\bibitem{Pham2009Continuous}
Huy\^{e}n Pham.
\newblock {\em Continuous-time stochastic control and optimization with financial applications}, volume~61 of {\em Stochastic Modelling and Applied Probability}.
\newblock Springer-Verlag, Berlin, 2009.

\bibitem{Raissi2019Physics}
M.~Raissi, P.~Perdikaris, and G.~E. Karniadakis.
\newblock Physics-informed neural networks: a deep learning framework for solving forward and inverse problems involving nonlinear partial differential equations.
\newblock {\em J. Comput. Phys.}, 378:686--707, 2019.

\bibitem{Ruthotto20}
Lars Ruthotto and Eldad Haber.
\newblock Deep neural networks motivated by partial differential equations.
\newblock {\em J. Math. Imaging Vision}, 62(3):352--364, 2020.

\bibitem{Ruthotto2020Machine}
Lars Ruthotto, Stanley~J. Osher, Wuchen Li, Levon Nurbekyan, and Samy~Wu Fung.
\newblock A machine learning framework for solving high-dimensional mean field game and mean field control problems.
\newblock {\em Proc. Natl. Acad. Sci. USA}, 117(17):9183--9193, 2020.

\bibitem{shi2024stochastic}
Zekun Shi, Zheyuan Hu, Min Lin, and Kenji Kawaguchi.
\newblock Stochastic taylor derivative estimator: Efficient amortization for arbitrary differential operators.
\newblock In {\em The Thirty-eighth Annual Conference on Neural Information Processing Systems}, 2024.

\bibitem{Sirignano2018DGM}
Justin Sirignano and Konstantinos Spiliopoulos.
\newblock D{GM}: a deep learning algorithm for solving partial differential equations.
\newblock {\em J. Comput. Phys.}, 375:1339--1364, 2018.

\bibitem{Sun2017Ito}
Yabing Sun, Jie Yang, and Weidong Zhao.
\newblock It\^o-{T}aylor schemes for solving mean-field stochastic differential equations.
\newblock {\em Numer. Math. Theory Methods Appl.}, 10(4):798--828, 2017.

\bibitem{Sun2018New}
Yabing Sun and Weidong Zhao.
\newblock New second-order schemes for forward backward stochastic differential equations.
\newblock {\em East Asian J. Appl. Math.}, 8(3):399--421, 2018.

\bibitem{Sun2020Explicit}
Yabing Sun and Weidong Zhao.
\newblock An explicit second-order numerical scheme for mean-field forward backward stochastic differential equations.
\newblock {\em Numer. Algorithms}, 84(1):253--283, 2020.

\bibitem{Sun2018Explicit}
Yabing Sun, Weidong Zhao, and Tao Zhou.
\newblock Explicit theta-schemes for mean-field backward stochastic differential equations.
\newblock {\em SIAM J. Numer. Anal.}, 56(4):2672--2697, 2018.

\bibitem{Sutton2018Reinforcement}
Richard~S. Sutton and Andrew~G. Barto.
\newblock {\em Reinforcement learning. {An} introduction}.
\newblock Adapt. Comput. Mach. Learn. Cambridge, MA: MIT Press, 2nd expanded and updated edition edition, 2018.

\bibitem{Tang2023DAS}
Kejun Tang, Xiaoliang Wan, and Chao Yang.
\newblock Das-pinns: A deep adaptive sampling method for solving high-dimensional partial differential equations.
\newblock {\em Journal of Computational Physics}, 476:111868, 2023.

\bibitem{Xiaoliang}
Xiaoliang Wan, Tao Zhou, and Yuancheng Zhou.
\newblock Adaptive importance sampling for deep {R}itz.
\newblock {\em Commun. Appl. Math. Comput.}, 7(3):929--953, 2025.

\bibitem{zhang2025shotgun}
Wenjun Xu and Wenzhong Zhang.
\newblock A deep shotgun method for solving high-dimensional parabolic partial differential equations.
\newblock {\em J. Sci. Comput.}, 104(2):69, 2025.

\bibitem{Yong13}
Jiongmin Yong.
\newblock Linear-quadratic optimal control problems for mean-field stochastic differential equations.
\newblock {\em SIAM J. Control Optim.}, 51(4):2809--2838, 2013.

\bibitem{Yong1999Stochastic}
Jiongmin Yong and Xun~Yu Zhou.
\newblock {\em Stochastic controls}, volume~43 of {\em Applications of Mathematics (New York)}.
\newblock Springer-Verlag, New York, 1999.
\newblock Hamiltonian systems and HJB equations.

\bibitem{Zang2020Weak}
Yaohua Zang, Gang Bao, Xiaojing Ye, and Haomin Zhou.
\newblock Weak adversarial networks for high-dimensional partial differential equations.
\newblock {\em J. Comput. Phys.}, 411:109409, 14, 2020.

\bibitem{Zhang2022FBSDE}
Wenzhong Zhang and Wei Cai.
\newblock F{BSDE} based neural network algorithms for high-dimensional quasilinear parabolic {PDE}s.
\newblock {\em J. Comput. Phys.}, 470:Paper No. 111557, 14, 2022.

\bibitem{Zhou2021Actor}
Mo~Zhou, Jiequn Han, and Jianfeng Lu.
\newblock Actor-critic method for high dimensional static {H}amilton-{J}acobi-{B}ellman partial differential equations based on neural networks.
\newblock {\em SIAM J. Sci. Comput.}, 43(6):A4043--A4066, 2021.

\bibitem{Zhou2024Solving}
Mo~Zhou and Jianfeng Lu.
\newblock Solving {Time}-{Continuous} {Stochastic} {Optimal} {Control} {Problems}: {Algorithm} {Design} and {Convergence} {Analysis} of {Actor}-{Critic} {Flow}.
\newblock Preprint, {arXiv}:2402.17208 [math.{OC}] (2024), 2024.

\bibitem{Zhou2025Policy}
Mo~Zhou and Jianfeng Lu.
\newblock A policy gradient framework for stochastic optimal control problems with global convergence guarantee.
\newblock {\em SIAM J. Control Optim.}, 63(4):2605--2631, 2025.

\bibitem{Zhu2025Optimal}
Yuhua Zhu, Yuming Zhang, and Haoyu Zhang.
\newblock Optimal-{PhiBE}: {A} {PDE}-based {Model}-free framework for {Continuous}-time {Reinforcement} {Learning}.
\newblock Preprint, {arXiv}:2506.05208 [math.{OC}] (2025), 2025.

\end{thebibliography}
 
\end{document}